\documentclass[a4paper, 12pt]{article}

\usepackage[T1]{fontenc} 
\usepackage{pict2e}
\makeatletter

\newcommand*\@KP@Large@frame[2]{%
	\setlength\unitlength{\fontdimen 22 #1\tw@}%
	\vrule \@width\z@ \@height 4\unitlength \@depth\tw@\unitlength
	\begin{picture}(6,2)(-3,-1)%
		\def\@KP@Radius     {3}%
		\def\@KP@Hole@radius{.5}%
		\def\@KP@Diameter   {6}%
		#2%
	\end{picture}%
}
\newcommand*\@KP@Small@frame[2]{%
	\setlength\unitlength{\fontdimen 22 #1\tw@}%
	\vrule \@width\z@ \@height \thr@@\unitlength \@depth\@ne\unitlength
	\begin{picture}(4,2)(-2,-1)%
		\def\@KP@Radius     {2}%
		\def\@KP@Hole@radius{.5}%
		\def\@KP@Diameter   {4}%
		#2%
	\end{picture}%
}

\newcommand*\@KP@Radius     {}
\newcommand*\@KP@Hole@radius{}
\newcommand*\@KP@Diameter   {}

\newcommand*\@KP@Shape@A{%
	\put(0,0){\circle{\@KP@Diameter}}%
}
\newcommand*\@KP@Shape@B{%
	\Line(-\@KP@Radius,\@KP@Radius )(\@KP@Radius,-\@KP@Radius)%
	\Line(-\@KP@Radius,-\@KP@Radius)(-\@KP@Hole@radius,-\@KP@Hole@radius)%
	\Line(\@KP@Radius ,\@KP@Radius )(\@KP@Hole@radius ,\@KP@Hole@radius )
}
\newcommand*\@KP@Shape@C{
	\cbezier(-\@KP@Radius,\@KP@Radius )(0,0)(0,0)(\@KP@Radius,\@KP@Radius )
	\cbezier(-\@KP@Radius,-\@KP@Radius)(0,0)(0,0)(\@KP@Radius,-\@KP@Radius)
}
\newcommand*\@KP@Shape@D{
	\cbezier(-\@KP@Radius,-\@KP@Radius)(0,0)(0,0)(-\@KP@Radius,\@KP@Radius)
	\cbezier(\@KP@Radius ,-\@KP@Radius)(0,0)(0,0)(\@KP@Radius ,\@KP@Radius)
}

\newcommand*\@KP@Shape@E{
	\cbezier(-\@KP@Radius,-\@KP@Radius)(0,0)(0,0)(\@KP@Radius,\@KP@Radius)
	\cbezier(-\@KP@Radius,-\@KP@Radius)(0,0)(0,0)(\@KP@Radius,\@KP@Radius)
}

\newcommand*\@KP@Atomic@mathpalette[1]{
	\mathinner{%
		\mathchoice{
			\linethickness{.6\p@}%
			\@KP@Large@frame \textfont {#1}
		}{
			\linethickness{.4\p@}
			\@KP@Small@frame \textfont {#1}
		}{
			\linethickness{.3\p@}
			\@KP@Small@frame \scriptfont {#1}
		}{
			\linethickness{.2\p@}
			\@KP@Small@frame \scriptscriptfont {#1}
		}
	}
}

\newcommand*\@KP@Shape@F{
	
	\put(-2, -3){\line(0,2){6}} 
	\put(-1, -3){\line(0,2){6}}

	\put(0,0){\circle*{0.2}} 
	\put(0.5,0){\circle*{0.2}} 
	\put(1,0){\circle*{0.2}}

	\put(2, -3){\line(0,2){6}} 
}

\newcommand*\KPF{\@KP@Atomic@mathpalette \@KP@Shape@F}

\newcommand*\KPA{\@KP@Atomic@mathpalette \@KP@Shape@A}
\newcommand*\KPB{\@KP@Atomic@mathpalette \@KP@Shape@B}
\newcommand*\KPC{\@KP@Atomic@mathpalette \@KP@Shape@C}
\newcommand*\KPD{\@KP@Atomic@mathpalette \@KP@Shape@D}
\newcommand*\KPE{\@KP@Atomic@mathpalette \@KP@Shape@E}
\makeatother

\usepackage{amsthm,amsmath,amssymb}
\usepackage{mathrsfs}
\usepackage{amsfonts}
\usepackage{graphicx}
\usepackage{pdfpages}
\usepackage{float}
\usepackage{url}
\usepackage{tikz-cd}

\usepackage{diagbox} 
\usepackage[utf8]{inputenc} 
\usepackage[T2A]{fontenc} 
\usepackage[russian,english]{babel} 
\usepackage{geometry}
\usepackage{hyperref}
\hypersetup{
	hypertex=true,
	colorlinks=true,
	linkcolor=blue,
	anchorcolor=blue,
	citecolor=blue}
\usepackage{subfigure}

\theoremstyle{plain}
\newtheorem{definition}{Definition}[section]
\newtheorem{example}{Example}

\newtheorem{theorem}{Theorem}[section]

\newtheorem{proposition}{Proposition}[section]
\newtheorem{lemma}{Lemma}[section]
\newtheorem{remark}{Remark}[section]

\usepackage{authblk}
\newcommand\keywords[1]{\textbf{Keywords}: #1}
\newcommand\msc[1]{\textbf{Mathematics Subject Classification}: #1}

\emergencystretch=\maxdimen
\usepackage{xcolor}
\usepackage[colorinlistoftodos]{todonotes}

\title{Multisacle Jones Polynomial and Persistent Jones Polynomial for Knot Data Analysis}

\author[1,2]{Ruzhi Song}
\author[1]{Fengling Li\thanks{fenglingli@dlut.edu.cn}}
\author[2]{Jie Wu}
\author[1]{Fengchun Lei}
\author[3,4,5]{Guo-Wei Wei}
\affil[1]{School of Mathematical Sciences, Dalian University of Technology, Dalian 116024, China}
\affil[2]{Beijing Institute of Mathematical Sciences and Applications, Beijing 101408, China}
\affil[3]{Department of Mathematics, Michigan State University, MI 48824, USA;}
\affil[4]{Department of Biochemistry and Molecular Biology, Michigan State University, MI 48824, USA}
\affil[5]{Department of Electrical and Computer Engineering, Michigan State University, MI 48824, USA}

\date{}

\begin{document}

	\maketitle

\begin{abstract}
	Many structures in science, engineering, and art can be viewed as curves in 3-space.
	The entanglement of these curves plays a crucial role in determining the functionality and physical properties of materials. 
	Many concepts in knot theory provide theoretical tools to explore the complexity and entanglement of curves in 3-space. 
	However, classical knot theory primarily focuses on global topological properties and lacks the consideration of local structural information, which is critical in practical applications.  
	In this work, two localized models based on the Jones polynomial, namely the multiscale Jones polynomial and the persistent Jones polynomial, are proposed. 
	The stability of these models, especially the insensitivity of the multiscale and persistent Jones polynomial models to small perturbations in curve collections, is analyzed, thus ensuring their robustness for real-world applications.
\end{abstract}


\msc{57K14, 92C10}

\section{Introduction}

Knot theory, a branch of mathematics that focuses on the study of mathematical knots, is primarily concerned with classifying and analyzing knots based on their essential properties under ambient isotopy \cite{crowell2012introduction}.
This approach allows mathematicians to disregard the specific manner in which knots are embedded in 3-space, emphasizing instead the invariants that remain unchanged under continuous deformations. 
There are various knot invariants, for example, the knot crossing number, the knot group \cite{crowell2012introduction}, Alexander polynomial \cite{Alexander1928Topological}, Jones polynomial \cite{Jones1985A}, knot Floer homology \cite{manolescu2016introduction}, and Khovanov homology \cite{khovanov2000categorification}.

Knot theory has applications across fields including physics \cite{ohtsuki2002quantum}, chemistry \cite{liang1994knots}, and biology \cite{sumners2020role,schlick2021knot,millett2013identifying}. 
In practical applications, however, two primary challenges arise: many structures do not form closed loops, and ambient isotopy can significantly alter local structures while preserving global knot characteristics. 
For example, open curves in 3-space, such as polymers \cite{qin2011counting, liu2018geometry}, textiles \cite{ricca2008topology}, chemical compounds \cite{panagiotou2019topological}, and biological molecules \cite{arsuaga2005dna, sulkowska2012conservation}, often exhibit local entanglement, which critically impacts their physical properties and functional roles. 
Topological invariants, functions that remain unchanged under ambient isotopy, are essential in analyzing knots and links \cite{freyd1985new,kauffman1990invariant,przytycki1987conway}.  However, these invariants do not extend to open curves, as open curves can be continuously deformed without requiring cutting or rejoining, making topological equivalence inapplicable.

In recent years, methods incorporating classical concepts from knot theory that are more applicable to practical problems have been proposed. 
Compared to topological data analysis (TDA), the concept of knot data analysis (KDA) was formally introduced in \cite{shen2024knot}.
Panagiotou and Plaxco \cite{panagiotou2020topological} demonstrated the utility of the Gauss link integral in protein entanglement, specifically for understanding protein folding kinetics and improving future folding models. 
Building on this, Baldwin and Panagiotou \cite{baldwin2021local} introduced a new measure of local topological and geometrical free energy based on the writhe and torsion of protein chains, highlighting its critical role in the rate-limiting steps of protein folding. 
Furthermore, Baldwin et al. \cite{baldwin2022local} extended these topological concepts to the study of the SARS-CoV-2 Spike protein, showing how local geometric features such as writhe and torsion influence its stability and behavior.
In a related effort, Shen et al. \cite{shen2024knot} introduced mGLI, a novel method leveraging the Gauss link integral to quantify the entanglement and topological complexity of both open and closed curves across various scales. 
This versatile approach has broad applications in analyzing curve structures in both physical and biological systems.

The Jones polynomial \cite{Jones1985A}, a fundamental invariant in knot theory, provides a polynomial measure of entanglement that distinguishes different types of knots by smoothing their crossings. Panagiotou and Kauffman \cite{panagiotou2020knot} extended this concept to an open curve, proposing a continuous measure of entanglement that converges to the classic Jones polynomial as the endpoints of an open curve approach each other. 
Barkataki and Panagiotou \cite{barkataki2022jones} further refined this by introducing the Jones polynomial for collections of curves, averaged across all projection directions. 
Building on these topological frameworks, Panagiotou and Kauffman \cite{panagiotou2021vassiliev} also applied Vassiliev invariants to quantify the complexity of both open and closed curves in 3-space. 
Additionally, Wang and Panagiotou \cite{wang2022protein} explored correlations between protein folding rates and topological measures, specifically, writhe, average crossing number (ACN), and the second Vassiliev invariant, to understand the behavior of native protein states.
In addition, Herschberg, Pifer, and Panagiotou \cite{herschberg2023computational} developed a computational tool called TEPPP, which quantifies topological complexity in systems like polymers, proteins, and periodic structures.

Utilizing the Jones polynomial framework for the collections of disjoint open or closed curves proposed by Barkataki and Panagiotou in \cite{barkataki2022jones}, this manuscript introduces two novel models: the multiscale Jones polynomial, represented by a characteristic matrix as detailed in Section \ref{subsectionMultiscale analysis of Jones polynomial}, and the persistent Jones polynomial, represented either through  weighted persistent barcodes or  weighted persistent diagrams, as described in Section \ref{subsection Persistent Jones Polynomial}. 
The weighted barcode was introduced by Cang and Wei in \cite{cang2020persistent}.
Both models have the capability to capture local and global entanglement properties in open or closed curve structures within 3-space and thus, effectively represent their topological characteristics. 
These models provide an enhanced approach for tackling problems in physical, biological, and chemical settings.

The stability of the two models is demonstrated. 
Specifically, small perturbations in the collection of disjoint open or closed curves lead to only minor changes in the characteristic matrix of the multiscale Jones polynomial and the weighted persistence diagram of the persistent Jones polynomial. 
More precisely, let \( L \) be a collection of disjoint open or closed curves in 3-space. 
Given a segmentation \( P_{n} = \{l_{1}, l_{2}, \dots, l_{n}\} \) for \( L \) with finite curve segments, where \( l_{i} \) represents a curve segment of \( L \), the segments \( l_{i} \) (\( 1 \leqslant i \leqslant n \)) can be connected end-to-end to reconstruct \( L \).
Consider a function \( f: L \to f(L) \), which induces a segmentation of \( f(L) \), denoted by \( f(P_{n}) = \{f(l_{1}), f(l_{2}), \dots, f(l_{n})\} \).

Suppose \( f: L \to f(L) \) is continuous and the difference \( \|f(L) - L\|_{\infty} < \varepsilon \) for sufficiently small \( \varepsilon > 0 \). 
There exist two characteristic matrices of the multiscale Jones polynomial for segmentations \( P_{n} \) and \( f(P_{n}) \), denoted by \( mJ(P_{n}) \) and \( mJ(f(P_{n})) \). 
The difference at each corresponding position between these two matrices remains sufficiently small. 
Additionally, the weighted Bottleneck distance between the weighted persistence diagrams of the persistent Jones polynomials for \( P_{n} \) and \( f(P_{n}) \) is also sufficiently small.

The proposed models are applied to B-factor prediction and the analysis of protein  $\alpha$-helix and $\beta$-sheet structures.
The B-factor, or Debye-Waller factor, is a critical metric in structural biology, representing the atomic displacement and flexibility within a protein structure, thus serving as an indicator of protein dynamics and stability. 
Traditional methods for predicting B-factors have had limitations in capturing topological information inherent in protein structures. 
To address this, we apply the multiscale Jones polynomial model to B-factor prediction, achieving prediction accuracies of \( 0.899 \), \( 0.808 \), and \( 0.720 \) for small, medium, and large protein sets \cite{park2013coarse}, respectively.
Our results on these three datasets  outperformed previous methods.
Additionally, the persistent Jones polynomial model is utilized to explore the structural properties of protein $\alpha$-helix and $\beta$-sheet segments, with visual representations provided through barcodes, highlighting the entanglement features across these secondary structures. 
The proposed multiscale Jones polynomial and persistent Jones polynomial models have potential for curve data analysis (CDA).

The article is organized as follows. 
In Section \ref{section The Jones polynomial of curves in 3-space}, the fundamental construction of the Jones polynomial for the collections of curves in 3-space is introduced.
In Section \ref{section Localization}, two new models of the Jones polynomial of curves in 3-space are established, namely the multiscale Jones polynomial (discussed in Section \ref{subsectionMultiscale analysis of Jones polynomial}) and the persistent Jones polynomial (described in Section \ref{subsection Persistent Jones Polynomial}).
In Section \ref{section Stability}, the stability of these two local models is demonstrated.
In Section \ref{section Application}, applications of the new models are presented, including B-factor prediction and the exploration of $\alpha$-helix and $\beta$-sheet structures.
In Section \ref{section Discussion}, a discussion on the selection of segmentation, localization, and stability issues is presented.

\section{The Jones polynomial of curves in 3-space}\label{section The Jones polynomial of curves in 3-space}

The Jones polynomial \cite{Jones1985A} is an important invariant in classical knot theory, recognized for its capacity to characterize the entanglement properties of knots and links. 
However, it is less applicable to practical scenarios often involving open curves in 3-space rather than closed curves. 
To address this limitation, Barkataki and Panagiotou \cite{barkataki2022jones} extended the concept to the collections of disjoint open or closed curves by defining a normalized version of the bracket polynomial, averaged over all projection directions. 
This adaptation enhances its relevance to real-world applications. 
Moreover, the Jones polynomial for collections of disjoint open or closed curves converges to the classical Jones polynomial as the endpoints of these curves approach one another. 
When projected onto a 2-dimensional plane, a collection of curves in 3-space forms a linkoid, which can be generalized to multi-component knotoids that describe open-ended knot diagrams. 
The theory of knotoids was initially introduced by Turaev \cite{turaev2012knotoids}, with further developments on linkoids elaborated in \cite{gugumcu2017new,gugumcu2017knotoids,gugumcu2019parity,manouras2021finite}.

\subsection{Segment Cycles}

Before defining the bracket polynomial of linkoids, it is essential to introduce the concept of segment cycles associated with a state.
Let $L$ be a linkoid diagram consisting of multiple components.
Let $G=\{1,2,3,\dots,2n\}$ denote the set of all endpoints (heads and legs) of $L$. 
A component of a linkoid with $n$ components is represented as $l_{2j - 1,2j}$, where $j\in\{1,2,\dots,n\}$.
The head-leg pairing forms a product of $n$ disjoint 2-cycles, denoted by $\hat{L}=(1,2)(3,4)\dots(2n - 1,2n)$.

Let $S$ be a state corresponding to a choice of smoothing over all crossing points in $L$. 
This state induces a pairing represented by the product of $n$ disjoint 2-cycles, $$\hat{S}=(s_1,s_2)(s_3,s_4)\dots(s_{2n - 1},s_{2n}),$$ where each $s_i\in G$ and each pair $(s_{2j - 1},s_{2j})$ for $j\in\{1,2,\dots,n\}$ represents the endpoints of a component in the state $S$.

For any endpoint $a\in G$, the set 
$$\mathrm{Orb}_S(a)=\{x\in G\mid x=(\hat{L}\circ\hat{S})^m(a),m\in\mathbb{Z}\}$$
is defined as the \textit{orbit} under the composition function $\hat{L}\circ\hat{S}$. 
The \textit{segment cycle} of an endpoint $a\in G$ is given by
$$\mathrm{Seg}(a)=\mathrm{Orb}_S(a)\sqcup\mathrm{Orb}_S(\hat{L}(a)).$$
It is notable that for any point $a\in G$, $\hat{L}(a)$ also belongs to the same segment cycle. 
Consequently, a segment cycle always contains an even number of elements.

\begin{lemma} \cite[Proposition 3.1]{barkataki2022jones}
	The number of segment cycles in a state $S$, denoted $|S|_{cyc}$, satisfies $1 \leqslant |S|_{cyc} \leqslant n$.
\end{lemma}

\begin{figure}[ht]
	\centering
	\begin{tikzpicture}

		\coordinate (A) at (0,2);
		\coordinate (HL_A) at (-1.5,1.35);
		\coordinate (J_A) at (1.5,1.3);
		\coordinate (HL_J_A) at (2.1,-0.1);
		\coordinate (J_HL_J_A) at (1.5,-1.5);
		\coordinate (HL_J_HL_J_A) at (0,-2.1);
		\coordinate (HL_J_A_) at (-2.08,0);

		\draw[thick] (A) arc[start angle=90, end angle=43, radius=2cm]; 
		\draw[thick, dashed] (A) arc[start angle=90, end angle=135, radius=2cm];
		
		\draw[thick, dashed] (J_A) arc[start angle=45, end angle=0, radius=2cm];
		
		\draw[thick] (HL_J_A) arc[start angle=0, end angle=-43, radius=2cm];
		
		\draw[thick, dashed] (J_HL_J_A) arc[start angle=-45, end angle=-90, radius=2cm];
		
		\draw[thick] (HL_A) arc[start angle=135, end angle=180, radius=2cm];

		\draw[thick] (HL_J_HL_J_A) arc[start angle=-90, end angle=-135, radius=2cm];

		\draw[thick, dotted] (HL_J_A_) arc[start angle=180, end angle=227, radius=2cm];

		\filldraw[black] (A) circle (2pt);
		\filldraw[black] (HL_A) circle (2pt);
		\filldraw[black] (J_A) circle (2pt);
		\filldraw[black] (HL_J_A) circle (2pt);
		\filldraw[black] (J_HL_J_A) circle (2pt);
		\filldraw[black] (HL_J_HL_J_A) circle (2pt);

		\node[above] at (A) {$a$};
		\node[left] at (HL_A) {$\hat{L}(a)$};
		\node[right] at (J_A) {$\hat{S}(a)$};
		\node[right] at (HL_J_A) {$\hat{L}(\hat{S}(a))$};
		\node[right] at (J_HL_J_A) {$\hat{S}(\hat{L}(\hat{S}(a)))$};
		\node[below] at (HL_J_HL_J_A) {$\hat{L}(\hat{S}(\hat{L}(\hat{S}(a))))$};
		
	\end{tikzpicture}
	
	\caption{Representation of the segment cycle of $a \in G$.}
	\label{segment_cycle_of_a}
\end{figure}
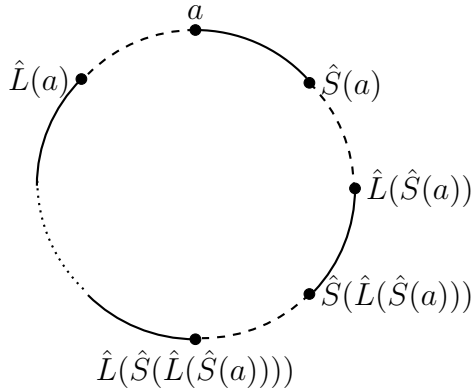

Consider a state $S$ of $L$ with the associated pairing $\hat{S}$. 
Let $\mathrm{Seg}(a)$ be a segment cycle in $S$ with $|\mathrm{Seg}(a)| = 2k$. 
This segment cycle can be represented by a circle marked with the $2k$ endpoints of $L$ (see Figure \ref{segment_cycle_of_a}). 
Let $a\in G$ be the initial point on the circle. 
The remaining $2k - 1$ endpoints are uniquely positioned in sequence on the circle as $\hat{S}(a)$, $\hat{L}(\hat{S}(a))$, $\hat{S}(\hat{L}(\hat{S}(a)))$, and so on, up to $\hat{L}(a)$.
It is important to note that the arcs connecting adjacent points in this circular representation alternate between the functions $\hat{S}$ and $\hat{L}$. 
Points connected by $\hat{S}$ belong to the same component in state $S$, while points connected by $\hat{L}$ belong to the same component in $L$.

\begin{figure}[htbp]
	\centering
	\subfigure[Hopf linkoid]
	{
		\begin{minipage}[b]{0.3\linewidth}
			\centering
			\begin{tikzpicture}

				\coordinate (P1) at (0,2);
				\coordinate (P2) at (0,0);
				\coordinate (P3) at (2,2);
				\coordinate (P4) at (2,0);

				\filldraw[black] (P1) circle (2pt) node[left] {1};
				\filldraw[black] (P2) circle (2pt) node[left] {2};
				\filldraw[black] (P3) circle (2pt) node[right] {3};
				\filldraw[black] (P4) circle (2pt) node[right] {4};

				\draw[-, thick] (P1) .. controls (1,1.5) and (1.6, 1.3) .. (1.2, 0.5);  
				\draw[-, thick] (1.0, 0.3) .. controls (0.6,0.0) and (0.4,0.0) .. (P2); 
				\draw[-, thick] (P3) .. controls (1.5,2) and (1.3,1.8) .. (1.1, 1.6);   
				\draw[-, thick] (0.9, 1.3) .. controls (0.7,1) and (0.3,0.6) .. (P4);   
				\draw[->,thick] (0.71,0.9) -- (0.71,0.85);
				\draw[->,thick] (1.32,1) -- (1.32,0.95);
			\end{tikzpicture}
		\end{minipage}
	}
	\subfigure[State $S_{1}$]
	{
		\begin{minipage}[b]{0.3\linewidth}
			\centering
			\begin{tikzpicture}
				
				\coordinate (P1) at (0,2);
				\coordinate (P2) at (0,0);
				\coordinate (P3) at (2,2);
				\coordinate (P4) at (2,0);

				\filldraw[black] (P1) circle (2pt) node[left] {1};
				\filldraw[black] (P2) circle (2pt) node[left] {2};
				\filldraw[black] (P3) circle (2pt) node[right] {3};
				\filldraw[black] (P4) circle (2pt) node[right] {4};

				\draw[-, thick] (P1) .. controls (0.7,1.5) and (1.4, 1.5) .. (P3);  
				\draw[-, thick] (P2) .. controls (-0.5,1.5) and (2.5,1.5) .. (P4); 				
			\end{tikzpicture}
		\end{minipage}
	}
	\subfigure[State $S_{2}$]
	{
		\begin{minipage}[b]{.3\linewidth}
			\centering
			\begin{tikzpicture}
			
				\coordinate (P1) at (0,2);
				\coordinate (P2) at (0,0);
				\coordinate (P3) at (2,2);
				\coordinate (P4) at (2,0);

				\filldraw[black] (P1) circle (2pt) node[left] {1};
				\filldraw[black] (P2) circle (2pt) node[left] {2};
				\filldraw[black] (P3) circle (2pt) node[right] {3};
				\filldraw[black] (P4) circle (2pt) node[right] {4};

				\draw[-, thick] (P1) .. controls (0.7,1.4) and (0.7, 0.7) .. (P2);  
				\draw[-, thick] (P3) .. controls (1.4,1.4) and (1.4,0.7) .. (P4); 			
			\end{tikzpicture}
		\end{minipage}
	}

	\caption{Hopf linkoid and its states}
	\label{Hopf_linkoid_state}
\end{figure}
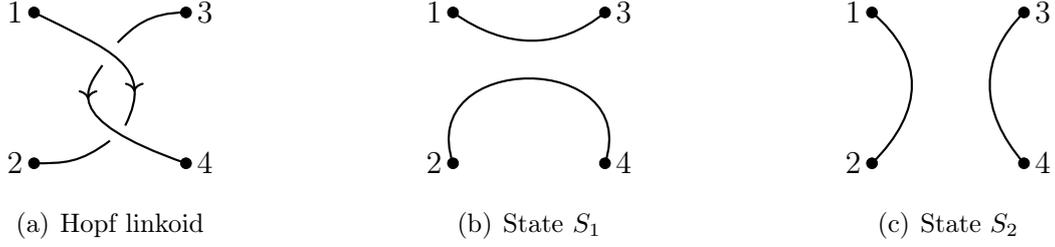

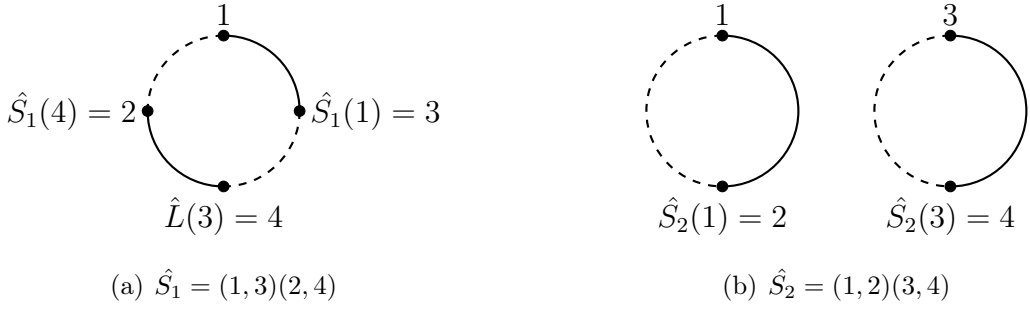
\begin{figure}[htbp]
	\centering
	
	\subfigure[$\hat{S_{1}} = (1,3)(2,4)$]
	{
		\begin{minipage}[b]{0.45\linewidth}
			\centering
			\begin{tikzpicture}

				\coordinate (A) at (0,1);
				\coordinate (HL_A) at (-1.5,1.35);
				\coordinate (J_A) at (1.5,1.3);
				\coordinate (HL_J_A) at (1,0);
				\coordinate (J_HL_J_A) at (1.5,-1.5);
				\coordinate (HL_J_HL_J_A) at (0,-1);
				\coordinate (HL_J_A_) at (-1,0);
			
				\draw[thick] (A) arc[start angle=90, end angle=0, radius=1cm];
				\draw[thick, dashed] (A) arc[start angle=90, end angle=180, radius=1cm];
				
				\draw[thick, dashed] (HL_J_A) arc[start angle=0, end angle=-90, radius=1cm];
				\draw[thick] (HL_J_A_) arc[start angle=180, end angle=270, radius=1cm];

				\filldraw[black] (A) circle (2pt);
				
				\filldraw[black] (HL_J_A) circle (2pt);
				
				\filldraw[black] (HL_J_HL_J_A) circle (2pt);
				\filldraw[black] (HL_J_A_) circle (2pt);

				\node[above] at (A) {$1$};
				
				\node[right] at (HL_J_A) {$\hat{S_{1}}(1)=3$};
				\node[left] at (HL_J_A_) {$\hat{S_{1}}(4)=2$};
				\node[below] at (HL_J_HL_J_A) {$\hat{L}(3)= 4$};
				
			\end{tikzpicture}
		\end{minipage}
	}
	\subfigure[$\hat{S_{2}} = (1,2)(3,4)$]
	{
		\begin{minipage}[b]{.45\linewidth}
			\centering
			\begin{tikzpicture}

				\coordinate (A) at (0,1);
				\coordinate (HL_A) at (-1.5,1.35);
				\coordinate (J_A) at (1.5,1.3);
				\coordinate (HL_J_A) at (1,0);
				\coordinate (J_HL_J_A) at (1.5,-1.5);
				\coordinate (HL_J_HL_J_A) at (0,-1);
				\coordinate (HL_J_A_) at (-1,0);
				\coordinate (A1) at (3,1);
				\coordinate (A2) at (3,-1);

				\draw[thick] (A) arc[start angle=90, end angle=-90, radius=1cm]; 
				\draw[thick, dashed] (A) arc[start angle=90, end angle=270, radius=1cm];
				
				\draw[thick] (A1) arc[start angle=90, end angle=-90, radius=1cm];
				\draw[thick, dashed] (A1) arc[start angle=90, end angle=270, radius=1cm];
				
				\filldraw[black] (A1) circle (2pt);
				
				\filldraw[black] (A2) circle (2pt);
				\node[above] at (A1) {$3$};
				
				\node[below] at (A2) {$\hat{S_{2}}(3)= 4$};

				\filldraw[black] (A) circle (2pt);
				
				\filldraw[black] (HL_J_HL_J_A) circle (2pt);
				
				\node[above] at (A) {$1$};
				
				\node[below] at (HL_J_HL_J_A) {$\hat{S_{2}}(1)= 2$};
				
			\end{tikzpicture}
		\end{minipage}
	}

	\caption{Segment cycles of two states}
	\label{Hopf_linkoid_segment_cycle}
\end{figure}

\begin{example}
	Consider the linkoid diagram $L$ as shown in Figure \ref{Hopf_linkoid_state}(a). 
	The set $G$ of all endpoints is $\{1,2,3,4\}$. 
	There are two states of $L$, $S_{1}$ and $S_{2}$, as shown in Figure \ref{Hopf_linkoid_state}(b) and (c). 
	The associated pairings $\hat{S_{1}}$ and $\hat{S_{2}}$ are represented by the permutations $(1,3)(2,4)$ and $(1,2)(3,4)$ respectively. 
	The segment cycles of states $S_{1}$ and $S_{2}$ are shown in Figure \ref{Hopf_linkoid_segment_cycle}. 
	Then $|S_{1}|_{cyc}=1$ and $|S_{2}|_{cyc}=2$.
\end{example}

\subsection{Jones polynomial}

The bracket polynomial of linkoids in $S^{2}$ or $\mathbb{R}^{2}$ is defined through an extension of the bracket polynomial of links.
The following initial conditions and diagrammatic relations are sufficient for the skein computation of the bracket polynomial of linkoids.

\begin{definition}\label{DefinitionOfLinkoid}
	Let $L$ be a linkoid diagram with $n$ components. 
	The bracket polynomial of the linkoid is uniquely determined by the following skein relation and initial conditions: 
	\[
	\left<L \cup \KPA\right> = (-A^{2} - A^{-2})\langle L \rangle,~
	\left< \KPB \right> = A \left< \KPC \right> + A^{-1} \left< \KPD \right>,~
	\left< \KPF \right> = (-A^{2} - A^{-2})^{|cyc|}
	\]
	where $|cyc|$ denotes the number of distinct segment cycles.
	
	The bracket polynomial of $L$ can be expressed as the following state sum expression:
	\[
	\left< L \right> = \sum_{S} A^{\sigma(S)} d^{|S|_{circ} - 1} d^{|S|_{cyc}},
	\]
	where $S$ is a state corresponding to a choice of smoothing over all crossing points in $L$; $\sigma(S)$ is the algebraic sum of the smoothing labels of $S$; $|S|_{circ}$ is the number of disjoint circles in $S$,  $|S|_{cyc}$ is the number of distinct segment cycles in $S$; and $d = (-A^{2} - A^{-2})$.
\end{definition}

The \textit{normalized bracket polynomial} is defined as follows:
\[
f_{L} = (-A^{-3})^{-Wr(L)} \left< L \right>,
\]
where $Wr(L)$ is the writhe of the linkoid diagram $L$.

Now consider curves in 3-space. 
A regular projection of curves fixed in 3-space can result in different linkoid diagrams depending on the projection direction chosen. 
Barkataki and Panagiotou in \cite{barkataki2022jones} define the bracket polynomial of curves in 3-space as the average of the bracket polynomial of a projection of the curve over all possible projection directions. 
This definition is made precise as follows:

\begin{definition}\cite[Definition 4.1.]{barkataki2022jones}
	Let $L$ denote a collection of disjoint open or closed curves in 3-space. 
	Let $(L)_{\xi}$ denote the projection of $L$ onto a plane with normal vector $\xi$. 
	The normalized bracket polynomial of $L$ is defined as follows:
	\[
	f_{L}=\dfrac{1}{4\pi}\int_{\xi\in S^{2}}(-A^{3})^{-Wr((L)_{\xi})}\langle (L)_{\xi}\rangle dS,
	\]
	where each $(L)_{\xi}$ is a linkoid diagram, and its bracket polynomial can be calculated using Definition \ref{DefinitionOfLinkoid}. 
	Note that the integral is taken over all vectors $\xi\in S^{2}$, excluding a set of measure zero (corresponding to the irregular projections). 
	This gives the Jones polynomial of a collection of disjoint open or closed curves in 3-space with the substitution $A = t^{-\frac{1}{4}}$.
\end{definition}

\begin{proposition}\cite[Proposition 4.1.]{barkataki2022jones}\label{PropertyJonesContinuous}
	\begin{itemize}
		\item For open curves, the Jones polynomial has real coefficients and is a continuous function of the curve coordinates.
		\item As the endpoints of the open curves in 3-space tend to coincide, the Jones polynomial tends to that of the corresponding link.
	\end{itemize}
\end{proposition}

\section{Localized Jones polynomials}\label{section Localization}

The Jones polynomial of a collection of disjoint open or closed curves in 3-space describes the entanglement of the curves within the entire collection. 
However, in many applications, it is desirable to extract the local structural information of the curves. 
Two methods for localizing the Jones polynomial are proposed to capture the entanglement of a collection of curves and meet the needs of practical applications.

\subsection{Multiscale Jones polynomial}\label{subsectionMultiscale analysis of Jones polynomial}

Let \( L \) be a collection of disjoint open or closed curves in 3-space. 
Given a segmentation \( P_{n} = \{l_{1}, l_{2}, \dots, l_{n}\} \) of \( L \), where each \( l_{i} \) represents a finite curve segment of \( L \), the segments \( l_{i} \), \( 1 \leq i \leq n \) can be sequentially connected to reconstruct \( L \).

To examine the entanglement properties between each curve segment and its neighboring segments, multiscale analysis is a suitable approach. 
A multiscale analysis requires defining a distance metric between the curve segments. 
The distance between segments, denoted \( d(l_{i}, l_{j}) \), can be specified by various metrics, depending on the application context. 
In this study, we define the distance $d(l_{i}, l_{j})$ as the upper bound of the Eulerian distances between a point of one curve segment and another curve segment for simplicity.

For any curve segment \( l_{i} \), consider the set of segments within \( P_{n} \) whose distances from \( l_{i} \) fall within the range \([r, R)\), where \( r \leqslant d(l_{i}, l_{j}) < R \). 
This set, which includes \( l_{i} \) itself, is denoted as
$$ P_{r,R}^{i} = \{l_{j} \in P_{n} ~|~ r \leqslant d(l_{i}, l_{j}) < R\} \cup \{l_{i}\}. $$

The Jones polynomial of the set of curve segments \( P_{r,R}^{i} \), denoted by \( JP_{r,R}^{i} \), quantifies the entanglement of the curve segment \( l_{i} \) with other segments as \( r \) and \( R \) vary. 
By selecting two sets of distance parameters, \( \{r_{1}, r_{2}, \dots, r_{m}\} \) and \( \{R_{1}, R_{2}, \dots, R_{m}\} \), with \( r_{i} < R_{i} \) for \( 1 \leq i \leq m \), we obtain a set of characteristic polynomials for \( l_{i} \):
\[
\{JP_{r_{1}, R_{1}}^{i}, JP_{r_{2}, R_{2}}^{i}, \dots, JP_{r_{m}, R_{m}}^{i}\}.
\]
Applying this procedure to all curve segments in \( P_{n} \), we obtain an \( n \times m \) matrix:
$$
\begin{pmatrix}
	JP_{r_{1}, R_{1}}^{1} & JP_{r_{2}, R_{2}}^{1} & \cdots & JP_{r_{m}, R_{m}}^{1}\\
	JP_{r_{1}, R_{1}}^{2} & JP_{r_{2}, R_{2}}^{2} & \cdots & JP_{r_{m}, R_{m}}^{2}\\
	\vdots & \vdots & \ddots & \vdots\\
	JP_{r_{1}, R_{1}}^{n} & JP_{r_{2}, R_{2}}^{n} & \cdots & JP_{r_{m}, R_{m}}^{n}\\
\end{pmatrix},
$$
capturing both local and global entanglement characteristics for the collection of disjoint open or closed curves \( L \).

Each matrix entry is a polynomial. 
For practical applications, the Jones polynomial can be evaluated at a specific parametrization, such as \( t = 10 \), resulting in an \( n \times m \) characteristic matrix for the segmentation \( P_{n} \) of \( L \):
$$
mJ(P_{n}) = 
\begin{pmatrix}
	JP_{r_{1}, R_{1}}^{1}(10) & JP_{r_{2}, R_{2}}^{1}(10) & \cdots & JP_{r_{m}, R_{m}}^{1}(10)\\
	JP_{r_{1}, R_{1}}^{2}(10) & JP_{r_{2}, R_{2}}^{2}(10) & \cdots & JP_{r_{m}, R_{m}}^{2}(10)\\
	\vdots & \vdots & \ddots & \vdots\\
	JP_{r_{1}, R_{1}}^{n}(10) & JP_{r_{2}, R_{2}}^{n}(10) & \cdots & JP_{r_{m}, R_{m}}^{n}(10)\\
\end{pmatrix},
$$
where each entry is a real number.

\begin{remark}
	For specific application contexts involving research objects that can be represented as a collection of disjoint open or closed curves \( L \), the choice of segmentation \( P_{n} = \{l_{1}, l_{2}, \dots, l_{n}\} \) plays a critical role. 
	Selecting an appropriate segmentation tailored to the application requirements allows for a more accurate capture of the entanglement properties inherent to the research objects, thus offering a more precise reflection of their characteristics. 
	Similarly, the choice of parameters \( \{r_{1}, r_{2}, \dots, r_{m}\} \) and \( \{R_{1}, R_{2}, \dots, R_{m}\} \) influences the robustness and accuracy of the final assessment of the objects' entanglement features.
\end{remark}

\subsection{Persistent Jones Polynomial}\label{subsection Persistent Jones Polynomial}

To address different application scenarios and to better capture the entanglement information of curves in 3-space, a second localization of the Jones polynomial has been proposed. 
This adaptation of the Jones polynomial for a collection of disjoint open or closed curves, aimed at quantifying entanglement complexity, serves as an extension of the classical Jones polynomial \cite{barkataki2022jones}. 

Let \( L \) be a collection of disjoint open or closed curves with a segmentation denoted by \( P_{n} = \{l_{1}, l_{2}, \dots, l_{n}\} \), where \( l_{i} \) is a segment of \( L \). The segments \( l_{i} \), \( 1 \leqslant i \leqslant n \), can be connected sequentially to reconstruct \( L \).
To effectively represent these multiple segments, both the \v{C}ech complex and the Vietoris-Rips complex are suitable methods, constructed here from the distance matrix:
$$
d(P_{n}) = 
\begin{pmatrix}
	0 & d(l_{1},l_{2}) & \cdots & d(l_{1},l_{n}) \\
	d(l_{2},l_{1}) & 0 & \cdots & d(l_{2},l_{n}) \\
	\vdots & \vdots & \ddots & \vdots \\
	d(l_{n},l_{1}) & d(l_{n},l_{2}) & \cdots & 0 \\
\end{pmatrix}
$$
of \( P_{n} = \{l_{1}, l_{2}, \dots, l_{n}\} \). 
Given the similarity between the \v{C}ech complex and the Vietoris-Rips complex, we focus on the Vietoris-Rips complex in the following discussion. 
Let \( r \) denote the variable parameter of the Vietoris-Rips complex.

\begin{definition}
	A critical value of the Vietoris-Rips complex is a real number \( r \) such that, for any sufficiently small \( \varepsilon > 0 \), the map \( K_{r-\varepsilon} \hookrightarrow K_{r+\varepsilon} \) is an inclusion but not an isomorphism, where \( K_{r} \) denotes the complex at Vietoris-Rips parameter \( r \).
\end{definition}

\begin{remark}
	For a segmentation \( P_{n} \) of \( L \) into finite curve segments, the Vietoris-Rips complex has a finite number of critical values.
\end{remark}

Let \( r_{0} < r_{1} < r_{2} < \dots < r_{m} \) represent the critical values of the Vietoris-Rips complex for a segmentation \( P_{n} \) of \( L \).
This generates a sequence of complexes from \( P_{n} \), forming a filtration \( \mathcal{F}(P_{n}) \):
\[
K_{r_0} \subsetneqq K_{r_1} \subsetneqq K_{r_2} \subsetneqq \dots \subsetneqq K_{r_m},
\]
where the final complex \( K_{r_{m}} \) is an \( (n-1) \)-simplex.
For any \( x < y \), let \( V_{x}^{y}: K_{x} \hookrightarrow K_{y} \) denote the inclusion map. 

\begin{lemma}\cite[Critical Value Lemma]{cohen2005stability}\label{Critical Value Lemma}
	If a closed interval \( [x, y] \) contains no critical value of the Vietoris-Rips complex, then \( V_{x}^{y}: K_{x} \to K_{y} \) is an isomorphism.
\end{lemma}

Within the filtration \( \mathcal{F}(P_{n}) \), consider a complex \( K_{r} \).
Each vertex \( v_{a} \subset K_{r} \) corresponds to a segment \( l_{a} \) from the curve segments in \( P_{n} \). 
An edge \( \{v_{a}, v_{b}\} \subset K_{r} \) indicates that the distance between segments \( l_{a} \) and \( l_{b} \) is less than the Vietoris-Rips parameter \( r \). 
A simplex \( \Delta = \{v_{a}, v_{b}, \ldots, v_{t}\} \subset K_{r} \) signifies that the pairwise distances among the corresponding segments \( \{l_{a}, l_{b}, \ldots, l_{t}\} \), collectively denoted by \( \Delta(P_{n}) \), are all less than \( r \).

\begin{definition}
	Let \( K \) be a simplicial complex. 
	The maximal faces of \( K \), with respect to inclusion, are named the \textit{facets} of \( K \). The simplicial complex \( K \), characterized by facets \( F_{1}, \ldots, F_{q} \), is denoted as
	$$ K = \langle F_{1}, \ldots, F_{q} \rangle, $$
	and the set \( \{F_{1}, \ldots, F_{q}\} \) is called the facet set of \( K \).
\end{definition}

Each complex within the filtration \( \mathcal{F}(P_{n}) \) can be described by its facets. 
Thus, the filtration \( \mathcal{F}(P_{n}) \) can be expressed as:
$$ \langle F^{r_{0}}_{1}, \ldots, F^{r_{0}}_{q_{0}} \rangle \subsetneqq \langle F^{r_{1}}_{1}, \ldots, F^{r_{1}}_{q_{1}} \rangle \subsetneqq \ldots \subsetneqq \langle F^{r_{m}}_{1}, \ldots, F^{r_{m}}_{q_{m}} \rangle. $$

\begin{definition}
	The \textit{birth} of a facet \( F \) in the filtration \( \mathcal{F}(P_{n}) \) is the smallest index \( r_{1} \) such that \( F \) appears as a facet in the complex \( K_{r_{1}} \) but not in \( K_{r_{1} - \varepsilon} \) for any sufficiently small \( \varepsilon > 0 \).
	
	The \textit{death} of a facet \( F \) in the filtration \( \mathcal{F}(P_{n}) \) is the largest index \( r_{2} \) such that \( F \) is a facet in \( K_{r_{2}} \) but not in \( K_{r_{2} + \varepsilon} \) for any sufficiently small \( \varepsilon > 0 \).
	
	The \textit{life-span} of a facet \( F \) in the filtration \( \mathcal{F}(P_{n}) \) is the interval \( [r_{1}, r_{2}] \).
\end{definition}

Therefore, \( \mathcal{F}(P_{n}) \) can be represented as a sequence of facets, each associated with a birth and death interval.
Similar to persistent barcodes in persistent homology, a barcode can represent the facets within a filtration. 
For any given dimension, each bar corresponds to a facet, with the bar’s start and end points indicating the birth and death of the related facet, respectively. 

\begin{definition}
	The barcode \( B(P_{n}) \) of the facets of the filtration \( \mathcal{F}(P_{n}) \) consists of horizontal line segments \( [r_{i}, r_{j}] \), where \( r_{i} \leq r_{j} \), representing the birth and death times of the associated facet.
\end{definition}

A barcode provides a visual representation of a filtration as a collection of horizontal line segments on a plane, where the horizontal axis corresponds to the parameter and the vertical axis represents an ordering of the facets.

Let \( \Delta = \{v_{a}, v_{b}, \ldots, v_{t}\} \) be a simplex in the filtration \( \mathcal{F}(P_{n}) \), and \( \Delta(P_{n}) = \{l_{a}, l_{b}, \ldots, l_{t}\} \) be the corresponding subset of segments in the segmentation \( P_{n} \).
The Jones polynomial of \( \Delta(P_{n}) \), denoted by \( J\Delta(P_{n}) \), can be regarded as a weight for the simplex \( \Delta \).
Thus, the Jones polynomial can be regarded as a weighting function over the filtration \( \mathcal{F}(P_{n}) \). 
The resulting weighted filtration is denoted by \( J\mathcal{F}(P_{n}) \).
We refer to the Jones polynomial weighted filtration \( J\mathcal{F}(P_{n}) \) as the \textit{persistent Jones polynomial} of segmentation \( P_{n} \) of \( L \).

Since the filtration \( \mathcal{F}(P_{n}) \) can be expressed by a barcode of facets, the persistent Jones polynomial of segmentation \( P_{n} \) of \( L \) can similarly be expressed by a barcode, with the Jones polynomials of the associated facets as weights.

The weights in the persistent Jones polynomial of \( P_{n} \) are polynomials. 
To enhance applicability in specific scenarios, setting the Jones polynomial variable \( t = 10 \) converts these weights into real numbers, producing a real-number weighted barcode \( BJ(P_{n})(10) \).

\section{Stability}\label{section Stability}

The stability of a model is characterized by the property that small perturbations in the collection of disjoint open or closed curves \( L \) result in only minor variations in the localized measures of the multiscale Jones polynomial and the persistent Jones polynomial.

Consider a continuous function \( f: L \to f(L) \), where \( f \) acts on a collection of disjoint open or closed curves \( L \) in 3-space. 
The difference between \( f(L) \) and \( L \) is measured using the supremum norm \( \|f(L) - L\|_{\infty} \), defined as:
\[
\|f(L) - L\|_{\infty} = \sup_{x \in L} |f(x) - x|.
\]

\begin{remark}\label{remark|J(L)(10)-J(f(L))(10)|varepsilon}
	The Jones polynomial evaluated at \( t = 10 \) can be considered a function on collections of curves in 3-space.  
	Let \( L \) be such a collection of curves in \( \mathbb{R}^3 \).  
	According to Proposition \ref{PropertyJonesContinuous}, let \( f : L \to f(L) \) be a continuous function. If \( \|f(L) - L\|_{\infty} < \varepsilon \) for all sufficiently small \( \varepsilon >0\), then \( |J(L)(10) - J(f(L))(10)| < \varepsilon_J \) for some sufficiently small \(\varepsilon_J > 0\).
\end{remark}

\subsection{Stability of Multiscale Jones Polynomial}

Let \( L \) be a collection of disjoint open or closed curves in 3-space, and let \( P_{n} = \{l_{1}, l_{2}, \dots, l_{n}\} \) denote a segmentation of \( L \) into \( n \) segments. Consider a continuous function \( f: L \to f(L) \) that induces a corresponding segmentation of \( f(L) \), represented by \( f(P_{n}) = \{f(l_{1}), f(l_{2}), \dots, f(l_{n})\} \).

\begin{proposition}\label{f(P_{r,R}^{i}) = f(P)_{r,R}^{i}}
	Suppose \( f: L \to f(L) \) is a continuous function such that \( \|f(L) - L\|_{\infty} < \varepsilon\) for all sufficiently small \( \varepsilon > 0 \). 
	Then, the two sets of curve segments \( f(P_{r,R}^{i}) \) and \( f(P)_{r,R}^{i} \) are equal:
	\[ f(P_{r,R}^{i})= f(P)_{r,R}^{i},\]
	where:
	\begin{itemize}
		\item \( f(P_{r,R}^{i}) \) is the image of the set \( P_{r,R}^{i} \) under the function \( f: L \to f(L) \),
		\item \( f(P)_{r,R}^{i} \) represents the set of curve segments in the segmentation \( f(P_n) \) such that their distance from \( f(l_i) \) is within \( [r, R) \), inclusive of \( f(l_i) \) itself.
	\end{itemize}
\end{proposition}

\begin{proof}
	By definition, we have \[ P_{r,R}^{i} = \{l_j \in P_n \mid r \leq d(l_i, l_j) < R \} \cup \{l_i\}, \]
	\[ f(P_{r,R}^{i}) = \{f(l_j) \in P_n \mid r \leq d(l_i, l_j) < R \} \cup \{f(l_i)\}, \] 
	and
	\[ f(P)_{r,R}^{i} = \{f(l_j) \in f(P_n) \mid r \leq d(f(l_i), f(l_j)) < R \} \cup \{f(l_i)\}. \]
	For any \( l_k \in P_{r,R}^{i} \), it holds that \( r \leq d(l_i, l_k) < R \).
	Since \( \|f(L) - L\|_{\infty} < \varepsilon \), the difference between each curve segment and its image under \( f \) is less than \( \varepsilon \):
	\( \| f(l_{i}) - l_{i} \|_{\infty} < \varepsilon \) and \( \| f(l_{k}) - l_{k} \|_{\infty} < \varepsilon \),
	and thus, the distances satisfy \( d(f(l_{i}), l_{i}) < \varepsilon \) and \( d(f(l_{k}), l_{k}) < \varepsilon \).
	
	\[
	\begin{tikzpicture}
		\centering
		\draw (0,0) node {$l_i$};
		\draw (2,0) node [above] {$d(l_i,l_k)$};
		\draw [dashed,<->] (0.8,0) -- (3.2,0);
		
		\draw (4,0) node {$l_k$};
		\draw [dashed,<->] (0,-0.6)--(0,-1.4);
		\draw (0,-2) node {$f(l_i)$};
		\draw [dashed,<->] (0.8,-2) -- (3.2,-2);
		\draw (4.3,-2) node {$f(l_k)$};
		\draw [dashed,<->] (4,-0.6)--(4,-1.4);
		\draw (0,-1) node [left] {$<\varepsilon$};
		\draw (4,-1) node [right] {$<\varepsilon$};
	\end{tikzpicture}
	\]
	
	Then, we have:
	\[
	d(l_i, l_k) - 2\varepsilon < d(f(l_i), f(l_k)) < d(l_i, l_k) + 2\varepsilon.
	\]
	Therefore,
	\[
	r - 2\varepsilon < d(f(l_i), f(l_k)) < R + 2\varepsilon,
	\]
	which implies \( f(l_k) \in f(P)_{r - 2\varepsilon, R + 2\varepsilon}^{i} \). Thus, \( f(P_{r,R}^{i}) \subseteq f(P)_{r - 2\varepsilon, R + 2\varepsilon}^{i} \).
	
	Similarly, we can show that:
	\[
	f(P)_{r - 2\varepsilon, R + 2\varepsilon}^{i} \subseteq f(P_{r - 4\varepsilon, R + 4\varepsilon}^{i}).
	\]
	Consequently, we obtain:
	\[
	f(P_{r,R}^{i}) \subseteq f(P)_{r - 2\varepsilon, R + 2\varepsilon}^{i} \subseteq f(P_{r - 4\varepsilon, R + 4\varepsilon}^{i}).
	\]
	Since \( \varepsilon > 0 \) is sufficiently small, we conclude that \( P_{r - 4\varepsilon, R + 4\varepsilon}^{i} = P_{r,R}^{i} \), \( f(P)_{r - 2\varepsilon, R + 2\varepsilon}^{i} = f(P)_{r,R}^{i} \), and
	\(
	f(P_{r,R}^{i}) = f(P)_{r,R}^{i}
	\).
\end{proof}

\begin{theorem}\label{theoremmultiscalanalysisJonespolynomial}
	Suppose \( f: L \to f(L) \) is a continuous function such that \( \|f(L) - L\|_{\infty} < \varepsilon \) for all sufficiently small \( \varepsilon > 0 \). 
	Consider two sets of distances \( \{r_1, r_2, \dots, r_m\} \) and \( \{R_1, R_2, \dots, R_m\} \). There exist two characteristic matrices for the segmentation \( P_n \) of \( L \) and for the segmentation \( f(P_n) \) of \( f(L) \), given by
	\[
	mJ(P_{n}) = 
	\begin{pmatrix}
		JP_{r_1, R_1}^{1}(10) & JP_{r_2, R_2}^{1}(10) & \cdots & JP_{r_m, R_m}^{1}(10) \\
		JP_{r_1, R_1}^{2}(10) & JP_{r_2, R_2}^{2}(10) & \cdots & JP_{r_m, R_m}^{2}(10) \\
		\vdots & \vdots & \ddots & \vdots \\
		JP_{r_1, R_1}^{n}(10) & JP_{r_2, R_2}^{n}(10) & \cdots & JP_{r_m, R_m}^{n}(10) \\
	\end{pmatrix},
	\]
	\[
	mJ(f(P_{n})) = 
	\begin{pmatrix}
		Jf(P)_{r_1, R_1}^{1}(10) & Jf(P)_{r_2, R_2}^{1}(10) & \cdots & Jf(P)_{r_m, R_m}^{1}(10) \\
		Jf(P)_{r_1, R_1}^{2}(10) & Jf(P)_{r_2, R_2}^{2}(10) & \cdots & Jf(P)_{r_m, R_m}^{2}(10) \\
		\vdots & \vdots & \ddots & \vdots \\
		Jf(P)_{r_1, R_1}^{n}(10) & Jf(P)_{r_2, R_2}^{n}(10) & \cdots & Jf(P)_{r_m, R_m}^{n}(10) \\
	\end{pmatrix}.
	\]
	
	Then, the difference between the corresponding entries in these two matrices is less than \( \varepsilon_{J} \), where \( \varepsilon_{J} \) is sufficiently small.
\end{theorem}

\begin{proof}
	We have
	$$\|P_{r_{j},R_{j}}^{i} - f(P_{r_j, R_j}^{i})\|_{\infty} \leq \|L - f(L)\|_{\infty} < \varepsilon.$$
	From Remark \ref{remark|J(L)(10)-J(f(L))(10)|varepsilon}, it follows that
	\[
	|JP_{r_j, R_j}^{i}(10) - Jf(P_{r_j, R_j}^{i})(10)| < \varepsilon_{J},
	\]
	where \( \varepsilon_{J} \) is sufficiently small.
	By Proposition \ref{f(P_{r,R}^{i}) = f(P)_{r,R}^{i}}, we know that \( f(P_{r,R}^{i}) = f(P)_{r,R}^{i} \).
	Therefore,
	\[
	|JP_{r_j, R_j}^{i}(10) - Jf(P)_{r_j, R_j}^{i}(10)| < \varepsilon_{J},
	\]
	where \( \varepsilon_{J} \) is sufficiently small.
\end{proof}

\begin{remark}
	The stability of the method used in \cite{shen2024knot} can be proved in a manner similar to that of Theorem \ref{theoremmultiscalanalysisJonespolynomial}. 
	Suppose \( f: L \to f(L) \) is a continuous function such that \( \|f(L) - L\|_{\infty} < \varepsilon \) for all sufficiently small \( \varepsilon > 0 \).
	
	For the Gauss linking integral, there exist two Gauss linking integral matrices for the segmentation \( P_{n} = \{l_{1}, l_{2}, \dots, l_{n}\} \) of a collection of disjoint open or closed curves \( L \) in 3-space, as well as for \( f(P_{n}) \) of \( f(L) \). These matrices are given by:
	\begin{align*}
		GL(P_{n}) &= 
		\begin{pmatrix}
			g(l_{1},l_{1}) & g(l_{1},l_{2}) & \cdots & g(l_{1},l_{n}) \\
			g(l_{2},l_{1}) & g(l_{2},l_{2}) & \cdots & g(l_{2},l_{n}) \\
			\vdots & \vdots & \ddots & \vdots \\
			g(l_{n},l_{1}) & g(l_{n},l_{2}) & \cdots & g(l_{n},l_{n})
		\end{pmatrix},
	\end{align*}
	\begin{align*}
		GL(f(P_{n})) &= 
		\begin{pmatrix}
			g(f(l_{1}),f(l_{1})) & g(f(l_{1}),f(l_{2})) & \cdots & g(f(l_{1}),f(l_{n})) \\
			g(f(l_{2}),f(l_{1})) & g(f(l_{2}),f(l_{2})) & \cdots & g(f(l_{2}),f(l_{n})) \\
			\vdots & \vdots & \ddots & \vdots \\
			g(f(l_{n}),f(l_{1})) & g(f(l_{n}),f(l_{2})) & \cdots & g(f(l_{n}),f(l_{n}))
		\end{pmatrix},
	\end{align*}
	where
	\begin{align}
		g(l_{i},l_{j}) = 
		\begin{cases}
			GL(l_{i},l_{j}), & \text{if } l_{i} \cap l_{j} = \emptyset, \\
			0, & \text{otherwise}.
		\end{cases}
	\end{align}
	Here, \( GL(l_{i},l_{j}) \) denotes the Gauss linking integral of curve segments \( l_{i} \) and \( l_{j} \).
	
	As stated in \cite[page 3]{panagiotou2020knot}, we can treat the Gauss linking integral of a curve as a continuous function of the curve's coordinates. 
	Similar to Theorem \ref{theoremmultiscalanalysisJonespolynomial}, we also have 
	\[ |g(l_{i},l_{j}) - g(f(l_{i}),f(l_{j}))| < \varepsilon_{GL} \]
	for sufficiently small \( \varepsilon_{GL} \).
\end{remark}

\subsection{Stability of the Persistence Jones Polynomial}

In this section, we state and prove the stability of the persistent Jones polynomial, which asserts that small changes in the collection of disjoint open or closed curves \( L \) lead to only small changes in the persistent Jones polynomial. 
As discussed in Section \ref{subsection Persistent Jones Polynomial}, the persistent Jones polynomial is represented by a weighted Jones polynomial filtration, which can be expressed as a Jones polynomial weight barcode of facets. 
For a given barcode, there is a corresponding diagram, where each bar in the barcode can be mapped to a point in the diagram. 
The $x$-coordinate of this point represents the birth time of the corresponding bar, while the $y$-coordinate represents its death time.

Let \( L \) be a collection of disjoint open or closed curves in 3-space, and let \( P_{n} = \{l_{1}, l_{2}, \dots, l_{n}\} \) denote a segmentation of \( L \) into \( n \) segments. Consider a continuous function \( f: L \to f(L) \) that induces a corresponding segmentation of \( f(L) \), represented by \( f(P_{n}) = \{f(l_{1}), f(l_{2}), \dots, f(l_{n})\} \). 
Let \( r_1 < r_2 < \dots < r_m \) be the critical values of the Vietoris-Rips complex of \( P_n \). We denote an interleaved sequence \( (b_i)_{i=0,1,\dots,m} \) such that \( b_{i-1} < r_i < b_i \) for all \( i \). We set \( b_{-1} = r_0 = -\infty \) and \( b_{m+1} = r_{m+1} = +\infty \).

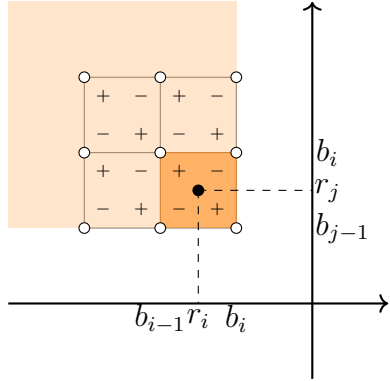
\begin{figure}[htbp]
	\centering
	\begin{tikzpicture}
		
		\draw[step=1cm, gray, very thin] (-3,1) grid (-1.0,3);

		\fill[orange, opacity=0.2] (-1, 1) rectangle (-4, 4);
		\fill[orange, opacity=0.5] (-1, 1) rectangle (-2, 2);

		\foreach \x in {-3, -2 , -1} {
			\foreach \y in {1, 2, 3} {
				\filldraw[fill=white, draw=black] (\x, \y) circle (2pt); 
			}
		}
		\filldraw[fill=black, draw=black] (-1.5, 1.5) circle (2pt);

		\node at (-2.75, 2.75) { \tiny$+$};
		\node at (-2.25, 2.75) { \tiny$-$};
		\node at (-1.75, 2.75) { \tiny$+$};
		\node at (-1.25, 2.75) { \tiny$-$};
		
		\node at (-2.75, 2.25) { \tiny$-$};
		\node at (-2.25, 2.25) { \tiny$+$};
		\node at (-1.75, 2.25) { \tiny$-$};
		\node at (-1.25, 2.25) { \tiny$+$};
		
		\node at (-2.75, 1.75) { \tiny$+$};
		\node at (-2.25, 1.75) { \tiny$-$};
		\node at (-1.75, 1.75) { \tiny$+$};
		\node at (-1.25, 1.75) { \tiny$-$};
		
		\node at (-2.75, 1.25) { \tiny$-$};
		\node at (-2.25, 1.25) { \tiny$+$};
		\node at (-1.75, 1.25) {\tiny$-$};
		\node at (-1.25, 1.25) {\tiny{$+$}};

		\draw[thick, ->] (-4,0) -- (1,0) node[anchor=north] {};
		\draw[thick, ->] (0,-1) -- (0,4) node[anchor=west] {};
		
		\draw[thin, dashed] (-1.5,1.5) -- (0,1.5);
		\draw[thin, dashed] (-1.5,1.5) -- (-1.5,0);

		\node at (-2, -0.2) {$b_{i-1}$};
		\node at (-1.5, -0.2) {$r_{i}$};
		\node at (-1, -0.2) {$b_{i}$};
		
		\node at (0.4, 1) {$b_{j-1}$};
		\node at (0.2, 1.5) {$r_{j}$};
		\node at (0.2, 2) {$b_{i}$};
	\end{tikzpicture}
	\caption{The multiplicity of the point \( (r_i, r_j) \) is the alternating sum at the corners of the lower right square. 
	When adding other multiplicities, cancellations between plus and minus signs occur.}
	\label{multiplicity}
\end{figure}

For two integers \( 0 \leqslant i < j \leqslant m+1 \) and a fixed integer \( k \), we define the \textit{multiplicity} of the pair \( (r_i, r_j) \) as
\[
\mu_i^j = \beta_{b_{i-1}}^{b_j} - \beta_{b_i}^{b_j} + \beta_{b_i}^{b_{j-1}} - \beta_{b_{i-1}}^{b_{j-1}},
\]
where \( \beta_x^y \) is the number of \( k \)-facets contained in \( K_x \) that remain in \( K_y \) for all \( -\infty \leqslant x \leqslant y \leqslant +\infty \).
To visualize this definition, consider \( \beta_x^y \) as the value of a function \( \beta \) at the point \( (r_{i}, r_{j}) \in \bar{\mathbb{R}}^2 \), where \( \bar{\mathbb{R}} = \mathbb{R} \cup \{-\infty, +\infty\} \). Thus, \( \mu_i^j \) is the alternating sum of \( \beta \) over the corners of the box \( [b_{i-1}, b_i] \times [b_{j-1}, b_j] \), as depicted in Figure \ref{multiplicity}.

Notice that if \( x \) and \( x' \) lie within the open interval \( (r_i, r_{i+1}) \), and \( y \) and \( y' \) lie within \( (r_{j-1}, r_j) \), then \( \beta_x^y = \beta_{x'}^{y'} \). Therefore, the multiplicities \( \mu_i^j \) are well-defined and always non-negative.

\begin{definition}
	The diagram \( D(P_n) \subset \bar{\mathbb{R}}^2 \) of the persistent Jones polynomial for \( P_n \) consists of points \( (r_i, r_j) \) with Jones polynomial weights, counted with multiplicity \( \mu_i^j \) for \( 0 \leq i < j \leq m+1 \), along with all points on the diagonal, which are counted with infinite multiplicity.
\end{definition}

Each off-diagonal point in the diagram represents the lifespan of a \( k \)-facet in the filtration. Similar to the weighted persistent barcode, the Jones polynomial corresponding to the set of curve segments for the \( k \)-facet can be used as the weight. This approach yields a weighted persistent diagram for the persistent Jones polynomial.

\subsubsection{Bottleneck Distance}

The \textit{Bottleneck distance} is a classical measure used to quantify the difference between two persistent diagrams. It naturally extends to the comparison of weighted persistent diagrams, allowing for a precise demonstration of the variations between two persistent Jones polynomials.

To better capture the differences between persistent Jones polynomials, we apply a slight modification to the traditional definition of the Bottleneck distance. 
Let \( \mathcal{C} \) and \( \mathcal{D} \) be two multi-sets of pairs \( (\langle a,b \rangle, w) \), where \( \langle a,b \rangle \) denotes an interval that can be any well-defined member of the set \( \{[a,b], [a,b), (a,b], (a,b)\} \), and \( w \in \mathbb{R} \) represents the weight of the interval \( \langle a,b \rangle \).

A \textit{matching} between the sets \( \mathcal{C} \) and \( \mathcal{D} \) is defined as a collection of pairs \( \chi = \{(I,J) \in \mathcal{C} \times \mathcal{D}\} \), where each element \( I \in \mathcal{C} \) and each element \( J \in \mathcal{D} \) appears in at most one pair within \( \chi \). 
A matching forms a bijection between a subset of \( \mathcal{C} \) and a subset of \( \mathcal{D} \). 
If a pair \( (I,J) \in \chi \), we say that \( I \) is \textit{matched} with \( J \). 
Conversely, if an element \( I \) does not appear in any pair, it is considered \textit{unmatched}.

The \textit{cost} \( c(I,J) \) of matching elements \( I = (\langle a,b \rangle, w_1) \) and \( J = (\langle c,d \rangle, w_2) \) is defined as follows:
\[
c(I,J) = \max\big\{|c-a|, |d-b|, |w_1 - w_2|\big\}.
\]
Similarly, the cost \( c(I) \) of leaving an unmatched element \( I \) is defined as:
\[
c(I) = \frac{b-a}{2}.
\]
Finally, the cost of a matching \( \chi \) is given by:
\[
c(\chi) = \max\left(\sup_{(I,J) \in \chi} c(I,J), \sup_{\text{unmatched } I \in \mathcal{C} \cup \mathcal{D}} c(I)\right).
\]

\begin{definition}
	The weighted Bottleneck distance between \( \mathcal{C} \) and \( \mathcal{D} \) is defined as
	\[
	d_B(\mathcal{C}, \mathcal{D}) = \inf\{c(\chi) \mid \chi \text{ is a matching between } \mathcal{C} \text{ and } \mathcal{D} \}.
	\]
\end{definition}

The modified Bottleneck distance increases the emphasis on weight factors compared to the classical Bottleneck distance.

\subsubsection{Stability}

Let \( L \) be a collection of disjoint open or closed curves in 3-space, and let \( P_{n} = \{l_{1}, l_{2}, \dots, l_{n}\} \) denote a segmentation of \( L \) into \( n \) segments. Consider a continuous function \( f: L \to f(L) \), that induces a segmentation of \( f(L) \), denoted by \( f(P_{n}) = \{f(l_{1}), f(l_{2}), \dots, f(l_{n})\} \).

There are two persistent Jones polynomials based on \( P_{n} \) and \( f(P_{n}) \), represented as \( J\mathcal{F}(P_{n}) \) and \( J\mathcal{F}(f(P_{n})) \), respectively. 
The weights of the facets in these persistent Jones polynomials are polynomials. 
By setting the Jones polynomial variable \( t = 10 \), the weight of each facet is converted into a real number. 
Thus, the converted persistent Jones polynomials can be denoted by \( J\mathcal{F}(P_{n})(10) \) and \( J\mathcal{F}(f(P_{n}))(10) \). 
These can be expressed using weighted persistent diagrams, denoted by \( D(P_{n}) \) and \( D(f(P_{n})) \).

Suppose \( \|f(L) - L\| < \varepsilon \) for all sufficiently small \( \varepsilon > 0 \). Then, for any \( l_{p_{i}} \in P_{n}(p) \), we have \( \|l_{p_i} - f(l_{p_i})\|_{\infty} < \varepsilon \). 
Let \( l_{p_{i}}, l_{p_{j}} \) be any two curve segments in \( P_{n}(p) \). 
Then, there exists
\[
\begin{tikzpicture}
	\centering
	\draw (0,0) node {$l_{p_{i}}$};
	\draw (2,0) node [above] {$d(l_{p_{i}},l_{p_{j}})$};
	\draw [dashed,<->] (0.8,0) -- (3.2,0);
	
	\draw (4,0) node {$l_{p_{j}}$};
	\draw [dashed,<->] (0,-0.6)--(0,-1.4);
	\draw (0,-2) node {$f(l_{p_{i}})$};
	\draw [dashed,<->] (0.8,-2) -- (3.2,-2);
	\draw (4.3,-2) node {$f(l_{p_{j}})$};
	\draw [dashed,<->] (4,-0.6)--(4,-1.4);
	\draw (0,-1) node [left] {$<\varepsilon$};
	\draw (4,-1) node [right] {$<\varepsilon$};
\end{tikzpicture}.
\]
Thus,
\[
d(l_{p_{i}}, l_{p_{j}}) - 2\varepsilon \leq d(f(l_{p_{i}}), f(l_{p_{j}})) \leq d(l_{p_{i}}, l_{p_{j}}) + 2\varepsilon.
\]
There is an important lemma, proved in \cite{cohen2005stability}.

\begin{lemma}\label{BoxLemma}\cite[Box Lemma]{cohen2005stability}
	For \( a < b < c < d \), let \( R = [a,b] \times [c,d] \) be a box in \( \mathbb{R}^{2} \), and let \( R_{2\varepsilon} = [a+2\varepsilon,b-2\varepsilon] \times [c+2\varepsilon, d-2\varepsilon] \) be the box obtained by shrinking \( R \) on all sides by \( 2\varepsilon \).
	It follows that:
	\[
	\#(D(P_{n}) \cap R_{2\varepsilon}) \leqslant \#(D(f(P_{n})) \cap R).
	\]
\end{lemma}

\begin{theorem}\label{BottleneckStability}
	Let \( L \) be a collection of disjoint open or closed curves, and let \( P_{n} = \{l_{1}, l_{2}, \dots, l_{n}\} \) denote a segmentation of \( L \) into \( n \) segments. Suppose \( f: L \to f(L) \) is a continuous function, such that \( \| f(L) - L \|_{\infty}< \varepsilon \) for all sufficiently small \( \varepsilon > 0 \). 
	Then, the weighted Bottleneck distance between the weighted persistent diagrams of persistent Jones polynomials, \( d_{B}(D(P_{n}), D(f(P_{n}))) \), is sufficiently small.
\end{theorem}

\begin{proof}	
	Let \( L \) be a collection of disjoint open or closed curves in 3-space, and let \( P_{n} = \{l_{1}, l_{2}, \dots, l_{n}\} \) denote a segmentation of \( L \) into \( n \) segments. Consider a continuous function \( f: L \to f(L) \), that induces a segmentation of \( f(L) \), denoted by \( f(P_{n}) = \{f(l_{1})\), \(f(l_{2})\), \(\dots, f(l_{n})\} \).
	There are two diagrams of the persistent Jones polynomial for \( P_{n} \) and \( f(P_{n}) \), denoted by \( D(P_{n}) \) and \( D(f(P_{n})) \).

	Consider the minimum distance between two distinct off-diagonal points or between an off-diagonal point and the diagonal:
	\[
	\delta_L = \min\{\|p - q\|_{\infty} \mid p \neq q \in D(P_{n}) - \Delta\}.
	\]
	Assuming \( \varepsilon > 0 \) is sufficiently small, we take \( \varepsilon < \delta_{L}/4 \).

	By drawing cubes of radius \( 2\varepsilon \) around points in \( D(P_{n}) \), we obtain a thickened diagonal plane along with a finite set of disjoint cubes, which are also disjoint from the thickened diagonal, as shown in Figure \ref{thickeneddiagonal}.

	\begin{figure}[htbp]
		\centering
		\begin{tikzpicture}
			
			\definecolor{lightorange}{RGB}{250,200,150}

			\filldraw[lightorange, draw = black] (0,1) -- (0,0) -- (1,0)--(6.3,5.3) -- (6.3,6.3) -- (5.3,6.3) -- cycle;
			
			\filldraw[fill=lightorange, draw=black ] (0.3, 0.3) rectangle (1.3, 1.3);
			\filldraw[fill=lightorange, draw=black ] (0.2, 0.2) rectangle (1.2, 1.2);
			\filldraw[fill=lightorange, draw=black ] (0.1, 0.1) rectangle (1.1, 1.1);
			\filldraw[fill=lightorange, draw=black ] (0, 0) rectangle (1, 1);
			\filldraw[fill=black, draw=black] (0.5, 0.5) circle (3pt);
			
			\filldraw[fill=lightorange, draw=black ] (5.3, 5.3) rectangle (6.3, 6.3);
			\filldraw[fill=lightorange, draw=black ] (5.2, 5.2) rectangle (6.2, 6.2);
			\filldraw[fill=lightorange, draw=black ] (5.1, 5.1) rectangle (6.1, 6.1);
			\filldraw[fill=lightorange, draw=black ] (5, 5) rectangle (6, 6);
			\filldraw[fill=black, draw=black] (5.5, 5.5) circle (3pt);

			\draw[dashed] (1.3,1.3) -- (5.5,5.5);

			\filldraw[fill=lightorange, draw=black ] (0, 5) rectangle (1, 6);
			\filldraw[fill=black, draw=black] (0.5, 5.5) circle (3pt);
			
			\filldraw[fill=lightorange, draw=black ] (0, 2.5) rectangle (1, 3.5);
			\filldraw[fill=black, draw=black] (0.5, 3) circle (3pt);
			
			\filldraw[fill=lightorange, draw=black ] (2.5, 5) rectangle (3.5, 6);
			\filldraw[fill=black, draw=black] (3, 5.5) circle (3pt);
		\end{tikzpicture}
		\caption{The shaded squares are centered at the black points of \( D(P_{n}) \).}
		\label{thickeneddiagonal}
	\end{figure}
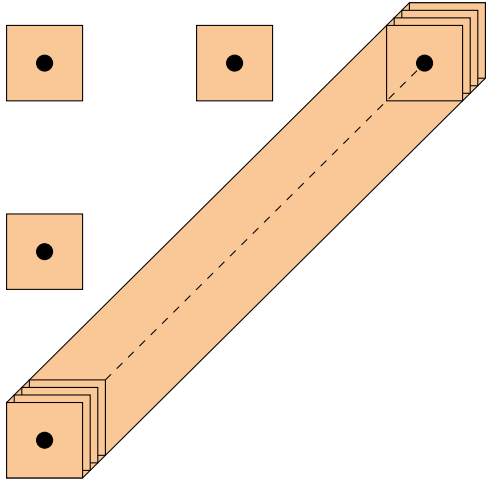

	Let \( \mu \) denote the multiplicity of a point \( p \) in \( D(P_{n}) \setminus \Delta \), and let \( \square_{2\varepsilon} \) represent the cube centered at \( p \) with radius \( 2\varepsilon \). According to Lemma \ref{BoxLemma},
	\[
	\mu \leq \#(D(f(P_{n})) \cap \square_{2\varepsilon}) \leq \#(D(P_{n}) \cap \square_{4\varepsilon}).
	\]
	Since \( 4\varepsilon < \delta_L \), \( p \) is the only point in \( D(P_{n}) \) within \( \square_{4\varepsilon} \), which implies that \( \#(D(f(P_{n})) \cap \square_{2\varepsilon}) = \mu \).

	Now let \( p = (x, y) \) be an off-diagonal point in the diagram \( D(P_{n}) \) with multiplicity \( \mu = 1 \). A corresponding collection of curve segments \( P_{n}(p) \) \(=\) \( \{l_{p_1}, l_{p_2}, \dots, l_{p_t}\} \) exists in \( P_{n} \). Applying the continuous function \( f \) to these segments yields the transformed collection \( f(P_{n}(p)) = \{f(l_{p_1}), f(l_{p_2}), \dots, f(l_{p_t})\} \).

	According to the definition of the persistent Jones polynomial, the point \( p = (x, y) \) indicates that \( P_{n}(p) = \{l_{p_1}, l_{p_2}, \dots, l_{p_t}\} \) forms a facet with birth \( x \) and death \( y \) in the filtration \( \mathcal{F}(P_{n}) \). This implies that any two segments in \( P_{n}(p) \) are at a distance less than \( x \), and any segment \( l \in P_{n} \setminus P_{n}(p) \) is at least a distance \( y \) from all \( l_{p_i} \in P_{n}(p) \), for \( 1 \leq i \leq t \).

	Since 
	\[
	d(l_{p_{i}}, l_{p_{j}}) - 2\varepsilon \leq d(f(l_{p_{i}}), f(l_{p_{j}})) \leq d(l_{p_{i}}, l_{p_{j}}) + 2\varepsilon,
	\]
	then any two curve segments in \( f(P_{n}(p)) \) are within a distance of \( x + 2\varepsilon \). Furthermore, for \( f(l) \in f(P_{n} \setminus P_{n}(p)) \) and \( f(l_{p_{i}}) \in f(P_{n}(p)) \), the distance is at least \( y - 2\varepsilon \).

	Thus, the segment collection \( f(P_{n}(p)) \) forms a facet within \( \mathcal{F}(f(P_{n})) \), corresponding to a point \( f(p) \) in \( D(f(P_{n})) \), within the region \( [x, x + 2\varepsilon] \times [y - 2\varepsilon, y] \), indicating that \( \|f(p) - p\|_{\infty} < 2\varepsilon \).
	Moreover, we have \( \#(D(f(P_{n})) \cap \square_{2\varepsilon}) = 1 \); hence, \( f(p) \) is the only point in \( D(f(P_{n})) \cap \square_{2\varepsilon} \).

	For off-diagonal points \( p^{1} = p^{2} = \dots = p^{\mu} = (x, y) \) with \( \mu > 1 \), this indicates \( \mu \) distinct facets in \( \mathcal{F}(P_{n}) \) with the same birth \( x \) and death \( y \), yet distinct curve segment collections \( P_{n}(p^1)\), \(P_{n}(p^2)\)\(, \dots, \)\(P_{n}(p^\mu) \). 
	Analogous to the case when \( \mu = 1 \), the images \( f(p^{1})\), \( f(p^{2})\)\(, \dots ,\)\(f(p^{\mu}) \) lie within \( [x, x + 2\varepsilon] \times [y - 2\varepsilon, y] \), yielding \( \|f(p^{i}) - p^{i}\|_{\infty} < 2\varepsilon \) for \( 1 \leq i \leq \mu \). 
	Moreover, we have that \( \#(D(f(P_{n})) \cap \square_{2\varepsilon}) = \mu, \) hence, \( f(p^{1}), f(p^{2}), \dots, f(p^{\mu}) \) are the \( \mu \) points in \( D(f(P_{n})) \cap \square_{2\varepsilon} \).

	The weights of \( p \) and \( f(p) \) are \( JP_{n}(p)(10) \) and \( Jf(P_{n}(p))(10) \), respectively. Given that 
	\[ \|f(P_{n}(p)) - P_{n}(p)\|_{\infty} \leq \|f(L) - L\|_{\infty} < \varepsilon, \]
	by Remark \ref{remark|J(L)(10)-J(f(L))(10)|varepsilon}, the weight difference satisfies \( |Jf(P_{n}(p))(10) - JP_{n}(p)(10)| < \varepsilon_J \), where \( \varepsilon_{J} > 0 \) is sufficiently small.

	After examining all off-diagonal points in \( D(P_{n}) \), the only points in \( D(f(P_{n})) \) that are not images \( f(p) \) for some \( p \in (D(P_{n})\setminus \Delta) \) lie beyond \( 2\varepsilon \) from \( D(P_{n}) \setminus \Delta \).

	Let \( q \in D(f(P_{n})) \) be a point for which there is no corresponding point \( p \in D(P_{n}) \) such that \( f(p) = q \). 
	Assume that the distance from \( q \) to \( \Delta \) is greater than \( 2\varepsilon \). Then, there exists a square \( \square^{q}_{2\varepsilon} \) centered at \( q \) with a radius of \( 2\varepsilon \) such that \( D(P_{n}) \cap \square^{q}_{2\varepsilon} = \emptyset \). This contradicts Lemma \ref{BoxLemma}, which states that \( 1 \leq \#(D(P_{n}) \cap \square^{q}_{2\varepsilon}) \neq 0 \). 
	Therefore, the distance from \( q \) to \( \Delta \) must be less than \( 2\varepsilon \).

	There exists a natural matching between \( D(P_{n}) \) and \( D(f(P_{n})) \), represented as \( \chi\)\( =\)\(\{(p, f(p)) \mid p \in D(P_{n}) \setminus \Delta \} \), with unmatched points in \( D(f(P_{n})) \) regarded as not corresponding to any point in \( D(P_{n}) \). Therefore, by the definition of the weighted Bottleneck distance,
	\[
	d_B(D(P_{n}), D(f(P_{n}))) < \max \{2\varepsilon, \varepsilon_J\},
	\]
	where \( \varepsilon \) and \( \varepsilon_J \) are sufficiently small, thus completing the proof.	
\end{proof}

In other words, the weighted persistent diagrams of persistent Jones polynomials are stable under small-amplitude or possibly irregular perturbations.

\section{Applications}\label{section Application}
\subsection{Multiscale Jones Polynomial for B-factor Prediction}

B-factors, also known as Debye-Waller factors, measure atomic displacement within protein structures, providing insight into molecular flexibility and stability. 
Analyzing B-factors enables a deeper understanding of protein dynamics and aids in predicting regions with high structural mobility, which is crucial for understanding protein function and interactions.

To eliminate the influence of non-relevant atomic information and to better capture the geometric and topological properties of the protein structure, each amino acid is represented by its \( C_{\alpha} \) atom. 
These \( C_{\alpha} \) atoms are sequentially connected to form a \( C_{\alpha} \) chain \( L \). 
Let \( C = \{c_0, c_1, \dots, c_n\} \) denote the set of \( C_{\alpha} \) atoms arranged in the sequence of the protein. 
The \( C_{\alpha} \) chain of the protein is considered a disjoint open curve. The segmentation of the \( C_{\alpha} \) chain is achieved by cutting at the midpoint between each \( C_{\alpha} \) atom and its adjacent \( C_{\alpha} \) atom, denoted by \( P_{n}=\{l_{0}, l_{1}, \dots , l_{n}\} \). The distance between two curve segments is defined as the distance between the \( C_{\alpha} \) atoms contained within the segments, that is, \( d(l_{i}, l_{j}) = d(c_{i}, c_{j}) \).

In this study, we select the radius $r$ to range from $4\mathring{\mathrm{A}}$ to $15\mathring{\mathrm{A}}$ with a step size of $0.25\mathring{\mathrm{A}}$. 
We set $R = (r + 1)\mathring{\mathrm{A}}$, so the radius $R$ ranges from $5\mathring{\mathrm{A}}$ to $16\mathring{\mathrm{A}}$. 
In total, the interception range is from $4\mathring{\mathrm{A}}$ to $16\mathring{\mathrm{A}}$.
Thus, there is a characteristic matrix for $P_{n}$ of the protein $C_{\alpha}$ chain $L$,
$$
\begin{pmatrix}
	JP_{4\mathring{\mathrm{A}}, 5\mathring{\mathrm{A}}}^{1}(10) & JP_{4.25\mathring{\mathrm{A}}, 5.25\mathring{\mathrm{A}}}^{1}(10) & \cdots & JP_{15\mathring{\mathrm{A}}, 16\mathring{\mathrm{A}}}^{1}(10)\\
	JP_{4\mathring{\mathrm{A}}, 5\mathring{\mathrm{A}}}^{2}(10) & JP_{4.25\mathring{\mathrm{A}}, 5.25\mathring{\mathrm{A}}}^{2}(10) & \cdots & JP_{15\mathring{\mathrm{A}}, 16\mathring{\mathrm{A}}}^{2}(10)\\
	\vdots & \vdots & \ddots & \vdots\\
	JP_{4\mathring{\mathrm{A}}, 5\mathring{\mathrm{A}}}^{n}(10) & JP_{4.25\mathring{\mathrm{A}}, 5.25\mathring{\mathrm{A}}}^{n}(10) & \cdots & JP_{15\mathring{\mathrm{A}}, 16\mathring{\mathrm{A}}}^{n}(10)\\
\end{pmatrix}.
$$
This choice is motivated by the fact that the average distance between \( C_{\alpha} \) atoms is approximately \( 3.8 \mathring{\mathrm{A}} \).
The selected radii scheme results in a powerful featurization method, providing rich representations of local protein structures.
To minimize the influence of overly complex machine learning models and highlight the effectiveness of the multiscale Jones Polynomial and avoid overfitting, we chose to use a Lasso regression model with parameter $0.16$ for B-factor prediction.

\begin{figure}[htbp]
	\centering
	\includegraphics[scale=0.4]{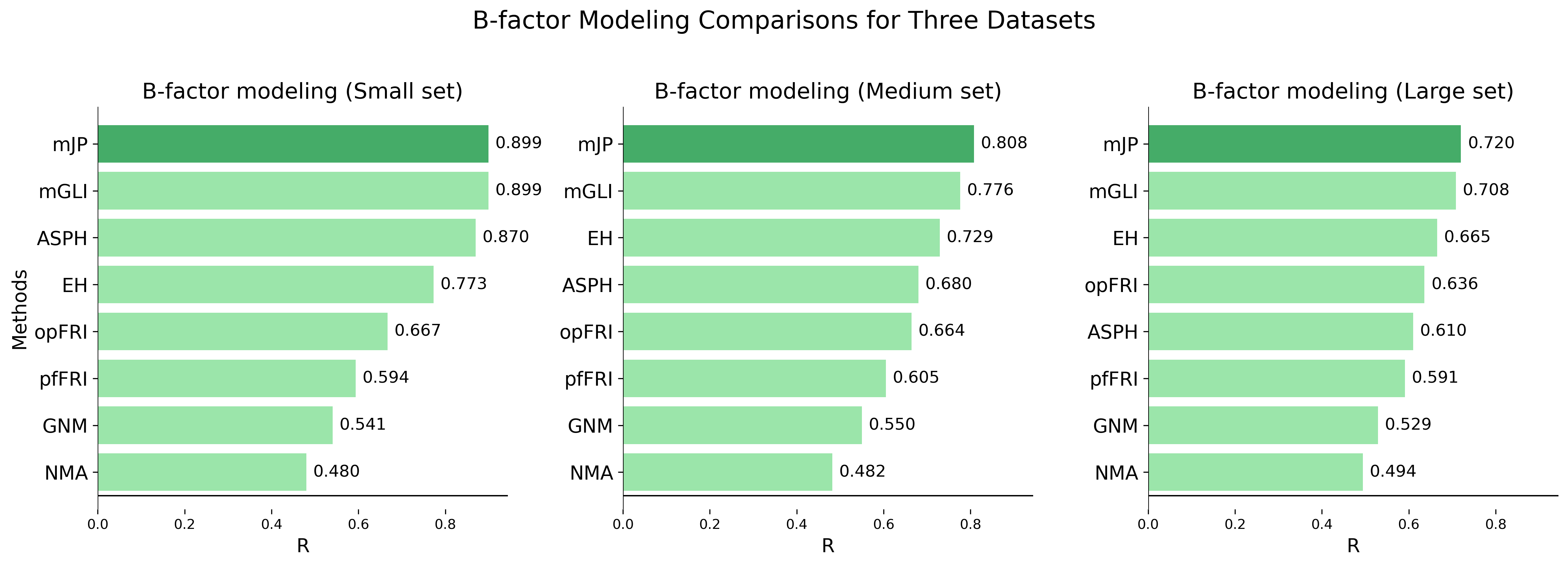}
	\caption{Comparison of B-factor predictions on three benchmark datasets between our Jones polynomial method and other approaches from the literature.}
	\label{B_factor_three_datasets}
\end{figure}

To validate the effectiveness of the multiscale Jones polynomial in predicting \( C_{\alpha} \) atom B-factors across proteins of varying sizes, we compared our method with several previous approaches, including mGLI \cite{shen2024knot}, EH \cite{cang2020evolutionary}, ASPH \cite{bramer2020atom}, opFRI \cite{opron2014fast}, pfFRI \cite{opron2014fast}, GNM \cite{park2013coarse}, and NMA \cite{park2013coarse}. The comparison was conducted on three protein sets from \cite{park2013coarse}, as shown in Figure \ref{B_factor_three_datasets}.  

The multiscale Jones polynomial method achieved average correlation coefficients of \( 0.899 \), \( 0.808 \), and \( 0.720 \) for small, medium, and large protein sets, respectively. 
Our results on these three datasets outperformed previous methods.

\begin{figure}[htbp]
	\centering
	\includegraphics[scale=0.5]{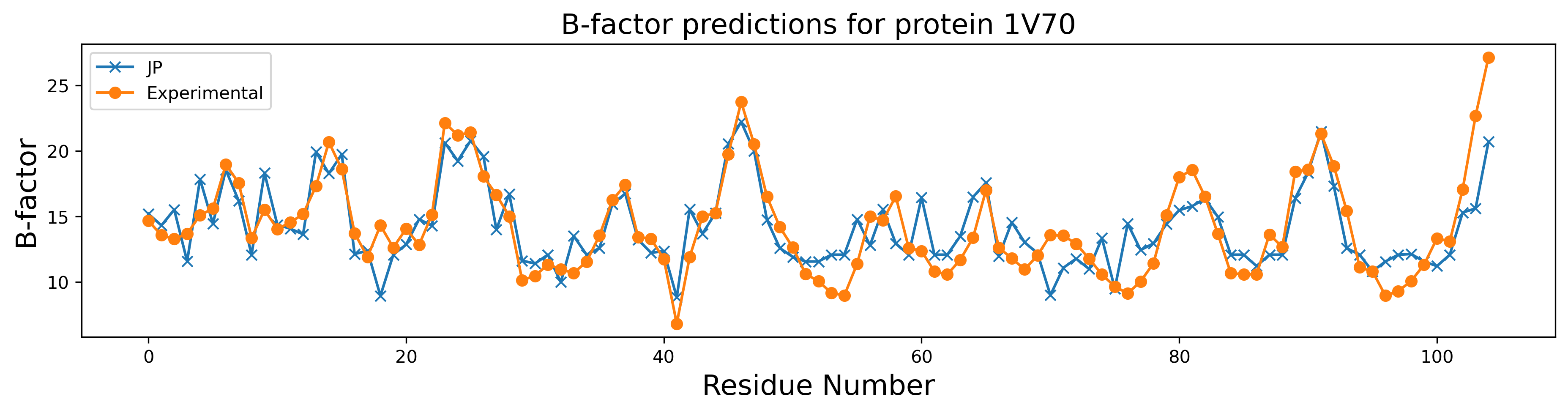}
	\caption{Comparison of B-factor values for protein $\mathrm{1V70}$, obtained from experiments and the Jones polynomial method. The $x$-axis represents the residue number, and the $y$-axis represents the B-factor value.}
	\label{B_facter_1v70}
\end{figure}

To further illustrate the performance of the multiscale Jones polynomial analysis, we present a case study of a potential antibiotic synthesis protein (PDBID: $\mathrm{1V70}$), which contains 105 residues. The B-factor values predicted by the multiscale Jones polynomial closely match the experimental values. A detailed comparison is provided in Figure \ref{B_facter_1v70}.  

Compared to traditional B-factor analysis methods, which focus on individual atoms, their spatial positions in 3-space, and the thermal motion and disorder of atoms within the protein structure, our approach effectively captures the torsional entanglement of the peptide chain at each \( C_{\alpha} \) atom's position by incorporating the Jones polynomial. 
This torsional entanglement of the peptide chain significantly influences the observed B-factor values.

The torsional entanglement of protein peptide chains, captured through both the multiscale and persistent Jones polynomials, provides critical insights into the structural and functional dynamics of proteins. By analyzing torsional entanglement, this approach reveals patterns of molecular flexibility and rigidity that allow for an in-depth understanding of protein stability and function. Furthermore, this method enhances our capacity to model and predict regions of structural mobility, offering potential applications for protein engineering and drug discovery.

\subsection{Barcodes of Persistent Jones Polynomial of $\alpha$-Helix and $\beta$-Sheet}\label{Barcodesofalphahelixandbetasheet}

In molecular biology, $\alpha$-helices and $\beta$-sheets are fundamental secondary structures in proteins, stabilized by hydrogen bonding patterns that contribute to the protein's overall stability and function. 
And $\alpha$-helices are typically more rigid than $\beta$-sheets. 
To explore the local structural complexity and stability of these structures, we employed topological analysis using the barcode representation of the persistent Jones polynomial. Using protein data from the Protein Data Bank (PDB), we demonstrated this approach with examples, including the analysis of an $\alpha$-helix chain consisting of 19 residues from the protein with PDB ID $\mathrm{1C26}$. 
Additionally, we extracted two parallel $\beta$-sheets from  protein 2JOX to explore their barcode representations of persistent Jones polynomials.

\begin{figure}[htbp]
	\centering

	\subfigure{
		\begin{minipage}[b]{.2\linewidth}
			\centering
			\includegraphics[scale=0.2]{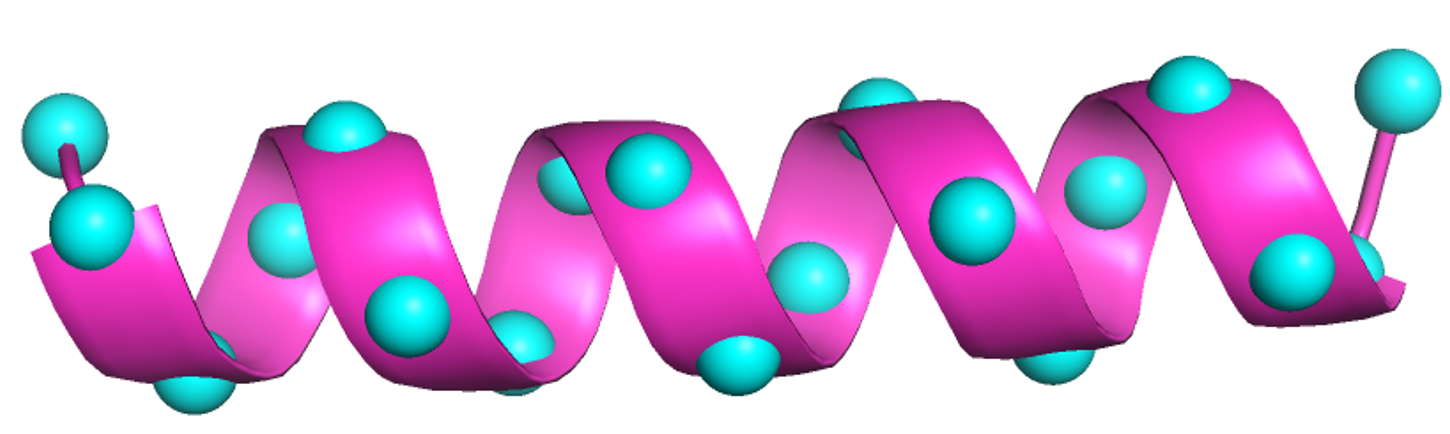}
		\end{minipage}
	}
	\subfigure{
		\begin{minipage}[b]{.75\linewidth}
			\centering
			\includegraphics[scale=0.3]{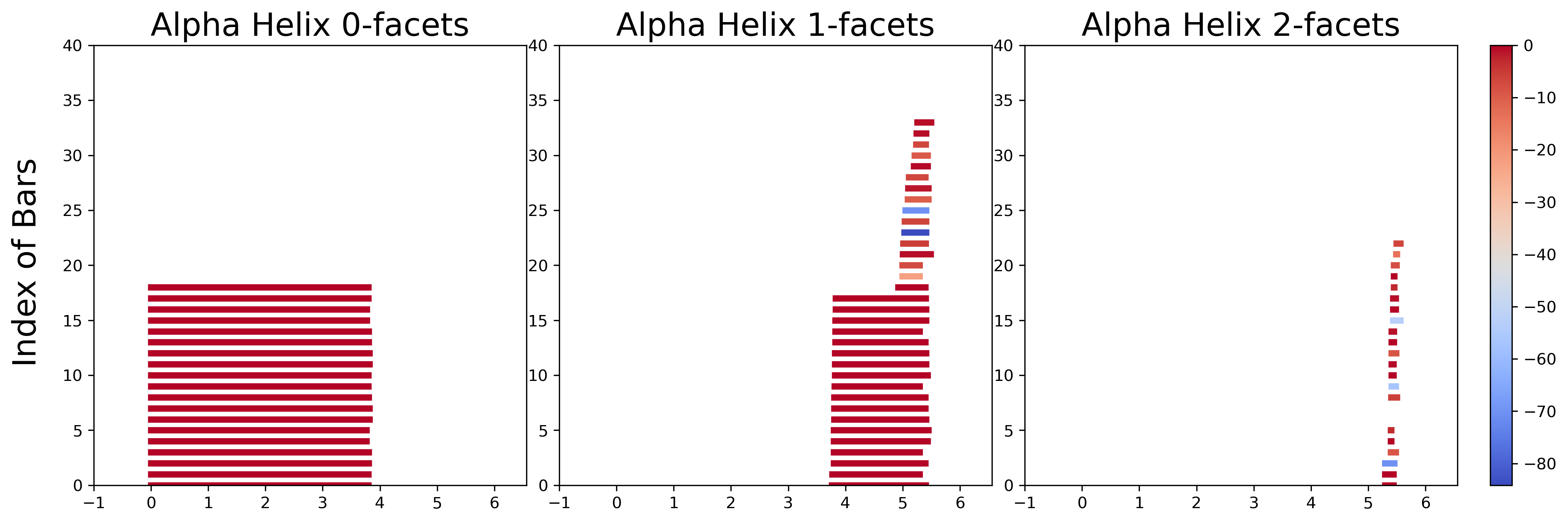}
		\end{minipage}
	}
	\caption{An alpha helix (left) and its persistent Jones polynomial barcodes (right). }
	\label{barcode_jones_alpha_helix}
\end{figure}

To eliminate the influence of non-relevant atomic information and better capture the geometric and topological properties of the helical structure, each amino acid is represented by its $C_{\alpha}$ atom, as shown in Figure \ref{barcode_jones_alpha_helix}. 
These $C_{\alpha}$ atoms are sequentially connected to form a $C_{\alpha}$ chain. 
The $C_{\alpha}$ chain of the $\alpha$-helix is considered a disjoint open curve. 
Segmentation of the $C_{\alpha}$ chain is achieved by cutting at the midpoint between each $C_{\alpha}$ atom and its adjacent $C_{\alpha}$ atom. 
The distance between two curve segments is defined as the distance between the $C_{\alpha}$ atoms contained within those segments.

\begin{figure}[htbp]
	\centering

	\subfigure{
		\begin{minipage}[b]{.2\linewidth}
			\centering
			\includegraphics[scale=0.2]{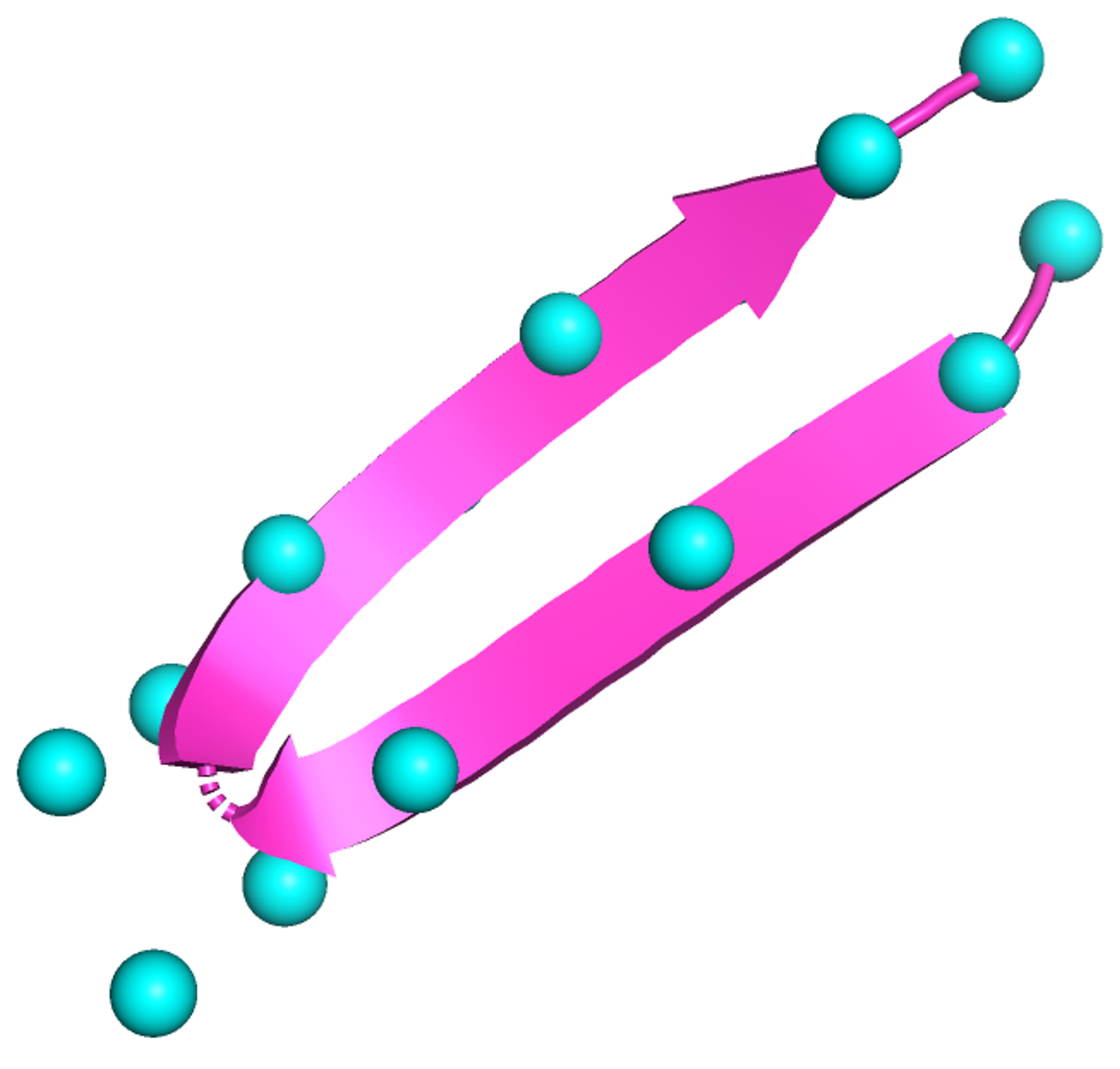}
		\end{minipage}
	}
	\subfigure{
		\begin{minipage}[b]{.75\linewidth}
			\centering
			\includegraphics[scale=0.3]{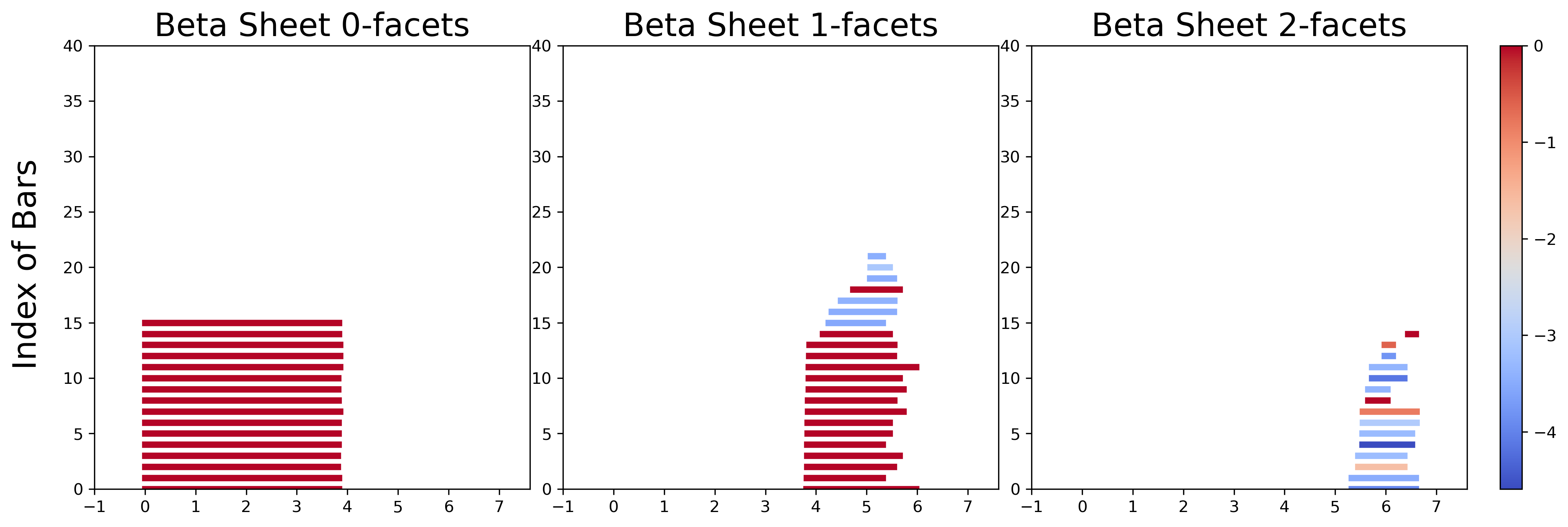}
		\end{minipage}
	}
	\caption{A beta sheet (left) and its persistent Jones polynomial barcodes (right).}
	\label{barcode_jones_beta_sheet}
\end{figure}

In the 0-facets panel of Figure \ref{barcode_jones_alpha_helix}, there are 19 bars with Jones polynomial weights of 0. 
Each bar has a length of approximately $3.8\mathring{\mathrm{A}}$, which represents the average distance between two $C_{\alpha}$ atoms. 
Additionally, in the 1-facets panel, there are 18 barcodes with similar birth times and persistent lengths. 
All of them begin around $3.8\mathring{\mathrm{A}}$ and persist until approximately $5.4\mathring{\mathrm{A}}$, each with a Jones polynomial weight of 0. 
These barcodes represent facets formed by two adjacent $C_{\alpha}$ atoms. 
Meanwhile, in the system of represented polylines, there is only one open polyline with two corner points. 
16 short-lived barcodes indicate facets formed by two non-adjacent $C_{\alpha}$ atoms. 
As shown in Figure \ref{barcode_jones_alpha_helix} (2-facets panel), the barcodes represent the persistence of facets formed by three $C_{\alpha}$ atoms and the entanglement of the corresponding polyline system.
Figure \ref{barcode_jones_beta_sheet} presents the persistent analysis using $C_{\alpha}$ atoms instead of amino acids.
Similar to the case of the $\alpha$-helix, the segmentation and its distance are defined in the same manner. 
There are 16 0-facets bars, indicating the presence of 16 $C_{\alpha}$ atoms. 
In the 1-facets panel, there are 14 facets formed by two adjacent $C_{\alpha}$ atoms and 8 facets formed by non-adjacent $C_{\alpha}$ atoms. 

The color of the barcodes reflects the value of the Jones polynomial weight, which indicates the difference of torsional entanglement of the set of curve segments in the system. 
It is important to note that the color gradient, whether tending towards red or blue, does not imply a higher or lower degree of entanglement complexity in the represented set of curve segments. 
Rather, a greater color difference between two bars suggests a greater difference between the sets of curve segments they represent.

From the colorbars in Figures \ref{barcode_jones_alpha_helix} and \ref{barcode_jones_beta_sheet}, it can be observed that the Jones polynomial weight for the $\alpha$-helix ranges from $-86$ to $0$, while for the $\beta$-sheet it ranges from $-6$ to $0$. 
This indicates a greater variability in the sets of curve segments represented by the facets within the $\alpha$-helix compared to the $\beta$-sheet, suggesting that the $\alpha$-helix exhibits more complex entanglement, as observed.

The persistent Jones polynomial effectively captures the torsional entanglement of secondary structures, such as $\alpha$-helices and $\beta$-sheets, within protein peptide chains. 
This torsional entanglement plays a critical role in analyzing protein structure and function, as it provides insights beyond atomic positions alone. 
By incorporating the Jones polynomial, our approach reveals important topological characteristics that are important to protein stability and functionality.

\section{Concluding remarks}\label{section Discussion}

In this study, the selection of segmentation \( P_{n} = \{l_{1}, l_{2}, \dots, l_{n}\} \) for a collection of curves \( L \) is crucial for capturing the topological and geometric characteristics of the curve structure \( L \). 
The segmentation serves as the basis for defining and calculating both the multiscale analysis of the Jones polynomial and the persistent Jones polynomial.
First, the outcomes of these two models depend not only on the spatial positions of the segments but also on their relative lengths in relation to the entire curve. 
When the segment length approaches zero, the results of the models tend toward triviality. 
Similarly, when the segment extends to cover the entire curve, the models recover global information. 
In both of these cases, the models are unable to extract meaningful  local information for spatial data. Thus, the choice of segmentation depends on the specific application.

Knot theory has traditionally focused on global invariants, but real-world applications often require local structural insights. 
Classical knot invariants mainly reflect global topology and fail to capture crucial local structural details in applications like molecular biology and highway crossing design. 
To address this gap, localized versions of invariants like the multiscale Jones polynomials and persistent Jones polynomials have been developed. 
These localized models decompose global invariants for analyzing local topology within the context of the entire structure, promoting the application of KDA or curve data analysis (CDA) in systems where both global and local structures matter.

Model stability is crucial for practical applications. 
In real-world data, noise and minor perturbations pose challenges. 
Stability ensures minor input changes do not cause disproportionate changes in calculated invariants. 
It is critical in biological or physical contexts. 
For multiscale Jones polynomial and persistent Jones polynomial models, we demonstrate stability under small perturbations. 
Minor adjustments in collection $L$ result in slight modifications to characteristic matrices and barcodes or diagrams. 
This stability makes the models reliable for structural topology and applicable in KDA or CDA for real-world data.  

The torsional entanglement of protein peptide chains, captured through both the multiscale and persistent Jones polynomials, provides crucial insights into the structure and function  of proteins. 
By analyzing torsional entanglement, this approach reveals patterns of molecular flexibility and rigidity that allow for an in-depth  understanding of protein structure and function. 
Furthermore, this method enhances our capacity to model and predict regions of structural mobility and reactivity, offering valuable implications for enzyme kinetics and protein engineering.

In Section \ref{Barcodesofalphahelixandbetasheet} of this manuscript, the barcodes of the persistent Jones polynomial represent the birth, death, and lifespan of facets in the complexes under filtration. 
In contrast, the barcodes in persistent homology capture the changes in homology classes of complexes during filtration, that is the changes of the generators of the homology groups of complexes (the number of generators corresponds to the Betti number).
Despite these differences, there are key similarities between the two concepts. Both rely on filtration and provide insights into data characteristics by examining the evolution of the complexes during filtration. 
Additionally, the bars in the barcodes of the persistent Jones polynomial represent facets formed by subsets of curve segments from the segmentation 
$P_{n} = \{l_{1}, l_{2}, \dots, l_{n}\}$ of the collection of curves $L$. 
Therefore, these bars are constructed based on the distance conditions between the curve segments in the segmentation $P_{n}$. 
Similarly, in the case of point cloud data, the barcodes of persistent homology are also constructed according to the distance conditions of the points in the point cloud data.

The multiscale Gauss link integral model \cite{shen2024knot} and the present  multiscale Jones polynomial and persistent Jones polynomial models represent solid advances in computational geometric topology. These approaches have great potential for real-world applications when they are paired with machine learning and artificial intelligence.

\section*{ACKNOWLEDGMENTS}
This work was supported in part by  grant
(No.12331003) of National Natural Science Foundation of China,  State Key
Laboratory of Structural Analysis, Optimization and CAE Software for Industrial
Equipment, and Beijing Institute of Mathematical Sciences and Applications.
The work of GWW was supported in part by NIH grants R01GM126189, R01AI164266, and R35GM148196, NSF grants DMS-2052983, DMS-1761320, DMS-2245903, and IIS-1900473, MSU Research Foundation, and Bristol-Myers Squibb 65109.

\bibliographystyle{plain}
\bibliography{refence}

\end{document}